\documentstyle[12pt]{amsart}
\newcommand{\bbc}{\mbox{$I \kern -6pt{C}$}}
\newcommand{\bbq}{\mbox{$I\kern -6pt{Q}$}}
\newcommand{\bbn}{\mbox{$I\kern -2pt{N}$}}
\newcommand{\bbz}{\mbox{$Z\kern -5pt{Z}$}}
\newcommand{\bbr}{\mbox{$I\kern -2pt{R}$}}
 
\newcommand{\catc}{\mbox{$\cal C$}}

\newcommand{\catcp}{\mbox{${\cal C}^p$}}
\newcommand{\catcpl}{\mbox{${\cal C}^+$}}

\newcommand{\catd}{\mbox{$\cal D$}}

\newcommand{\cate}{\mbox{$\cal E$}}
\newcommand{\catm}{\mbox{$\cal M$}}

\newcommand{\catlc}{\mbox{\bf LC}}
\newcommand{\catgen}{\mbox{\bf Gen}}
\newcommand{\cattop}{\mbox{\bf Top}}

\newcommand{\dom}[1]{\mbox{dom($#1$)}}
\newcommand{\domk}[1]{\mbox{dom\kern 2pt{$#1$}}}
\newcommand{\codk}[1]{\mbox{cod\kern 2pt{$#1$}}}

\newcommand{\morcp}[2]{\mbox{Mor$_{{\cal C}^{p}}(#1,#2)$}}

\newcounter{pfenumi}
\newenvironment{pflist}{\begin{list}
	{\rm (\arabic{equation}.\alph{pfenumi})}
	  {\usecounter{pfenumi}
	  \setlength{\labelsep}{.08 in}
	  \setlength{\rightmargin}{.15 in}
	  \setlength{\leftmargin}{.40 in} }
	\refstepcounter{equation} }{
	\end{list}  }
\newenvironment{onecond}{
	\begin{list}{\rm (\arabic{equation})}
	  {\usecounter{pfenumi}
	  \setlength{\labelsep}{.08 in}
	  \setlength{\rightmargin}{.15 in}
	  \setlength{\leftmargin}{.40 in} }
	\refstepcounter{equation} }{\end{list}  }

\newenvironment{pfproof}{\vspace{4mm} \noindent
	{\bf Proof.}~}{$\Box $\vspace{4mm}}

\newtheorem{theorem}{Theorem}
\newtheorem{proposition}[theorem]{Proposition}
\newtheorem{definition}[theorem]{Definition}
\newtheorem{lemma}[theorem]{Lemma}
\newtheorem{corollary}{Corollary}
\newtheorem{remark}{Remark}
\newtheorem{conjecture}[theorem]{Conjecture}

\begin{document}

\title{Existence of Orbifolds IV:  \\
 	Examples}

\author{Paul Feit}
\thanks{Much of this paper and its
predecessor was drafted during two summer visits to MSRI, Berkeley.
Research at MSRI is supported in part by
NSF Grant \# DMS-9022140} 
\address{University Of Texas at Permian Basin  \\
	Odessa, Texas  \\
	79762 }
\maketitle

\begin{abstract}
This work concludes a series of four papers on the
foundational theory of orbifolds and stacks.  We
apply the abstract theory, developed in its
predecessors, to orbifolds derived from manifolds.
Specifically, we show how the very concrete topological
base spaces associated to such orbifolds can be described and
manipulated in our universal language.  At the same time, we 
interpret our many categorical axioms in several explicit
contexts.
\end{abstract}

\section*{Introduction}
In a series of works, the author has developed a Existence Criterion
for categories of orbifolds.  In this concluding
paper, we offer
applications of the abstract machinery.  We begin with topological
spaces relevant to orbifolds, and develop enough point-set theory to
prove that the hypotheses of our universal propositions apply.  We also
give explicit meaning to their conclusions.

This paper owes much of its subject matter to Ms. Dorette Pronk, who,
as of this writing, is preparing for her doctorate.  The author's
original motivation arose from algebraic geometry.  Specifically, our
hope was to simplify and complete formulations for algebraic spaces
and stacks.  This starting point essentially forced a categorical
perspective upon us.  Ms. Pronk pointed out that there are many more
accessible, geometric theories of orbifolds.  She challenged us to
explain certain empirical observations and practical questions in
the context of our theory.  This paper is, primarily, a series of
answers to her questions.  Some comments are complete, some are partial.

The contents are roughly as follows.  We introduce a subcategory
\catgen~of the category of topological spaces.  To avoid problems with
the axioms of set theory, we limit it as follows:  for each set $S$,
we let $\catgen [S]$ be the subcategory of objects which have a
cover by open subsets of cardinality less than or equal to the
cardinality of $S$.  (When $S=\bbr$, the subcategory includes most
standard objects of study.)  This category is then examined from
several perspectives.
\begin{pflist}
\item	We assign to \catgen~a {\em pseudogeometric\/}
	topology, and verify many abstract conditions discussed
	in previous papers.  Actually, the majority of our
	universal constraints are translations of classical
	point-set ideas to the categorical level.  Each individual
	identity is simple to check.
\item	A universal theorem from \cite{O3} now assigns to \catgen~a
	pseudo\'etale topology.  That is,
	we assign to it a topology from which our plus construction
	creates orbifold-type objects.  For technical reasons, it is
	convenient to restrict the topology to $\catgen [S]$.
	(There are formal questions about \catgen~which we
	conjecture to be resolvable.)  The category
	$\catgen [S]^{+}$ will contain our first explicit, non-trivial
	orbifolds.  These include objects with mirrored boundaries and
	cone points, as discussed in \cite{WT}.
\item	Although members of $\catgen [S]$ have no differentiable or
	analytic structure, the explicit construction here has obvious
	analogues for a category of manifolds of any kind ($S=\bbr$).
	Consequently, our work shows to define pseudo\'etale topologies
	for differential or analytic manifolds.
\item	Inside $\catgen [S]$, we give examples of a group action
	$G$ on an object $M$ such that non-trivial members of $G$
	have fixed points but, in the categorical sense of \cite{O1},
	the quotient map
	$q:M\longrightarrow G\backslash M$ is Galois (ie., discrete,
	overlays absolutely and uniformly, etc.).
\item	Many classical orbifolds are regarded as some sort of
	a topological base space plus extra information.  Our theory
	does not use sheaves of structure to model objects.  We have
	an alternative perspective which allows for base spaces in
	contexts where they exist.  In our language, to say that each
	object in a category \catc~has a topological aspect is to
	assign to it a continuous functor
	$\Gamma :\catc \longrightarrow \catgen [S]$ which is
	{\em faithful.}  If $f$ and $g$ are distinct
	\catc -morphisms between two objects, then
	$\Gamma (f)$ and $\Gamma (g)$ are continuous functions
	between ``underlying'' base spaces; moreover, if
	$\Gamma (f)$ and $\Gamma (g)$ agree as functions,
	then $f=g$ in \catc .  We prove that, under an elementary
	hypothesis, every lift of $\Gamma$ (ie., $\Gamma ^{+}$ ,
	$(\Gamma ^{+})^{+}$ ,...) is also faithful.
\item	Let \catc~be a category which is being used to
	generate orbifolds (such as a category of manifolds of
	some type).  Let $M,N\in \catc$, and let $G$ and $H$ be
	groups with actions on $M$ and $N$, respectively.  Then
	$G\backslash M$ and $H\backslash N$ exist in \catcpl .
	What is a morphism
	$G\backslash M\longrightarrow H\backslash N$?  We
	consider a question of pathologies from \cite{GS}.
\end{pflist}
Throughout the paper, we use the notation, terminology
and theorems of earlier works.  We must assume that the
reader is familiar with, or has access to these papers.

Section~\ref{topol} defines the topological context \catgen .
(We emphasize certain subcategories, which will not be
discussed in this introduction.)  As hinted in previous works,
a vital concept is that a theory of orbifold-type objects should
begin with a topological category in which the class of
morphisms is restricted.  The definition of
{\em diffuse} continuous function requires some point-set topology.

We need topologies for \catgen .  Each object in \catgen~has
a topology in the point-set sense.  Our universal framework
discusses topology, but in a categorical sense.
Section~\ref{geometric} defines a
{\em pseudogeometric} topology for \catgen .  This amounts to
(1) rewriting point-set concepts, like subsets and connectedness,
in the abstract language and (2) verifying many categorical
restrictions.  Although numerous, each individual condition of
our universal axiom set translates to an easy statement in the
concrete situation.

The pseudogeometric topology is too elementary to generate
orbifolds.  In \cite{O3}, a method for deriving a pseudo\'etale
topology from a pseudogeometric topology was introduced.  It is
the latter construct which, in unison with the plus functor,
generates orbifold-type objects.  It is the pseudo\'etale topology
that occupies the our attention for the rest of the paper.

Definition of the derived pseudo\'etale is forbiddingly formal.
However, we have developed an infrastructure of theory for the
concept.  Section~\ref{pseudo} reviews what is known about
pseudo\'tale morphisms, and where there are gaps.  We conjecture
that some of the omissions can be resolved by further work.  Other
difficulties, such as the fact that \catgen~is not closed under
descent, are part of the theory.

We seek a context in which there is a finite group $G$ acting on
a manifold $M$ such that non-trivial members of $G$ have fixed
points and yet the quotient map
$q:M\longrightarrow G\backslash M$ has many of the properties
usually exhibited only by quotients for discrete actions.  Those
properties were discussed, at the level of categories, in \cite{O1}.
We characterized particularly good quotients as being
{\em Galois.}  Inside \catgen , we can actually exhibit such
group actions for the first time.  In this context, the challenge is
to find $G$ and $M$ for which the morphism $q$ is pseudo\'etale.
Theorem~\ref{xtheo2} gives an explicit hypothesis under which
this occurs.  It requires some effort to prove that theorem.
Section~\ref{prod} starts the proof with a study of fibered
products in \catgen .  Then, Section~\ref{fullquot} brings in
several universal lemmas established in earlier papers to
complete the argument.

Section~\ref{faith} changes topics.  It focuses on the following issue:
\begin{onecond}
\label{wff}
\item	Suppose $f$ and $g$ are morphisms between two
	differentiable or analytic orbifolds of some kind.  Assume that
	the underlying functions of $f$ and $g$ (that is, the continuous
	maps they determine between base spaces) agree.  Show that $f=g$.
\end{onecond}
We interpret (\ref{wff}) to state that the functor which sends an
orbifold to its underlying topological space is faithful.  We offer a
hypothesis under which faithfulness of a topological model is
preserved by the plus construction.

Section~\ref{whatis} interprets the abstract definition in an
explicit context.  Suppose $G$ and $H$ are groups which act,
respectively, on affine objects $M$ and $N$ such that
$G\backslash M$ and $H\backslash N$ exist as orbifold-type
objects.  There is a definition of morphism from
$G\backslash M\longrightarrow H\backslash N$ in our theory.
If $M$ and $N$ are both some type of manifold, then
there exist topological quotients
$G\backslash M$ and $H\backslash N$ which are quite explicit.
This leads to a deceptive principle: an orbifold morphism
$F:M\longrightarrow N$ factors to an orbifold morphism
$f:G\backslash M\longrightarrow H\backslash N$ if and only
if the underlying continuous function of $F$ factors through
the topological quotients.  Mathematicians have discovered that
this idea leads to contradictions.  We illustrate the falsity of the
principle by considering a pathology based on an example in
Schwarz's thesis \cite{GS}.

The author is especially grateful to MSRI in Berkeley.  Virtually
all of this paper, and much of its predecessor, were drafted
during visits to MSRI in the summers of 1993 and 1994.
The staff helped to arrange a synchronous stay by Dorette Pronk,
whose concerns motivate most of this paper.

\section{Topological Terminology}
\label{topol}
Let \cattop~be the category of topological spaces (in which the morphism
class consists of all continuous functions).  In the present paper,
we define a space $X\in \cattop$ to be
{\em locally connected\/} if each point has a basis consisting of connected
neighborhoods.  We do
{\em not\/} require that a locally connected space be Hausdorff.  Let
\catlc~denote the class of locally connected topological spaces.
Connected components are always closed; if $X\in \catlc$, then its
connected components are open as well.  Denote the closure of a subset
$C\subseteq X$ by $\overline{C}$.

Let $X\in \cattop$.  We say $X$ is
{\em locally Hausdorff\/}  if it has a cover by open subsets, each of which
is Hausdorff in the subset topology.  A class of locally Hausdorff objects
will appear later.

Let $X\in \cattop$.  By a
{\em subset element\/} of $X$, we mean a pair $(U,I)$ where $U$ is an open
subset of $X$ and $I$ is a subset of $U$.  We call an element $(U,I)$
{\em closed\/} if $I$ is closed in $U$.  Suppose
$f:X\longrightarrow Y$ is a continuous function between topological spaces.
For $(U,I)$ a subset element of $Y$, define the subset element
{\em $f^{-1}$(U,I),\/} or
{\em pullback of $(U,I)$ along $f$\/}, of $X$ to be $(f^{-1}U,f^{-1}I)$.

Now, we introduce some non-standard terminology.

\begin{definition}
Let $X\in \catlc$.  An open subset $U$ of $X$ is called
{\em Z-dense\/} (in $X$) if, for every non-empty open connected subset
$C$ of $X$, $U\cap C$ is non-empty and connected.  (One may regard
these as a generalization of sets which are dense in a Zariski topology,
as used in algebraic geometry.)  A subset $I\subseteq X$ is called
{\em negligible\/} (in $X$) if its complement is Z-dense.
\end{definition}

The following comments are tautological consequences of the definition.

\begin{proposition}
Let $X\in \catlc$.

(A)  A Z-dense subset is dense and open.

(B)  The intersection of a finite family of Z-dense open subsets
is a Z-dense open subset.

(C)  Let $U,V\subseteq X$ be open subsets such that $U\subseteq V$.
If $U$ is Z-dense, then
$V$ is Z-dense.

(D)  Let $U,V\subseteq X$ be open subsets such that $U\subseteq V$.  If
$V$ is Z-dense in $X$ and $U$ is Z-dense in $V$, then $U$ is Z-dense in $X$.

(E)  Let $U,V\subseteq X$ be open subsets.  If $U$ is Z-dense in $X$,
then $U\cap V$ is a Z-dense open subset of $V$.
\end{proposition}

\begin{corollary}
\label{lcor1}
Let $X\in \catlc$.

(A)  A negligible subset of $X$ is closed and nowhere dense.

(B)  The union of a finite family of negligible subsets of $X$ is negligible.

(C)  A closed subset of a negligible set is negligible.

(D)  Let $I$ be a negligible subset of $X$ and let $J$ be a negligible
subset of $X-I$.  Then $I\cup J$ is negligible in $X$.

(E)  Let $I,V\subseteq X$ be subsets.  If
$V$ is open and $I$ is negligible in $X$, then $I\cap V$ is a
negligible subset of $V$.
\end{corollary}

\begin{pfproof}
Part (C) of the Lemma requires the elementary fact that if
$D$ is a connected subset of $X\in \cattop$, and if
$D\subseteq E\subseteq \overline{D}$, then $E$ is connected.
All other claims are tautological.
\end{pfproof}

It is easy to find a local characterization of the property of being negligible.

\begin{lemma}
\label{ylem2}  Let $X\in \catlc$, and let $I\subseteq X$ be a
closed subset.  Then the following two conditions are equivalent.
\begin{pflist}
\label{zlist1}
\item	For each $x\in I$, there is a basis of neighborhoods
	${\cal V}_{x}$ at $x$ such that, for each $V\in {\cal V}_{x}$, both
	$V$ and $V-I$ are connected and non-empty.
\item	$I$ is negligible in $X$.
\end{pflist}
\end{lemma}

\begin{pfproof}
Obviously the problem is to prove (\ref{zlist1}.b) from
(\ref{zlist1}.a).  Hereafter, assume (\ref{zlist1}.a).

As a first step, observe that, trivially, $I$ has no interior.
We finish the proof by contradiction.

Let $C$ be a connected, non-empty open subset of $X$ for which
$C-I$ is not both connected and non-empty.  By our first step,
$C-I$ must be a non-empty, disconnected set.  Then there are two
non-empty closed subsets $A$ and $B$ of $C-I$ such that
$A\cap B=\emptyset$ and $A\cup B=C-I$.  Let $A^{\ast}$ and $B^{\ast}$
be the respective closures of $A$ and $B$
{\em with respect to $C$.}  Then $A^{\ast} \cup B^{\ast} =C$, because
$I$ is nowhere dense.  Since $C$ is connected, there must be
$x\in A^{\ast} \cap B^{\ast} $.  Clearly, $x\in I$.  Let
$V\in {\cal V}_{x}$ such that $V\subseteq C$.  Then $V-I$ is connected.
and non-empty.  Consequently, $V-I$ intersects only one of the sets
$A,B$.  Without loss of generality, assume
$(V-I)\cap A=\emptyset$.  Now  $V\cap A=(V-I)\cap A=\emptyset$.
But, since $x$ is in the closure of $A$, the latter statement is impossible.
\end{pfproof}

\begin{definition}
Let $X\in \cattop$.  A subset element $(U,I)$ is called
{\em negligible\/} if, in the subset topology, $U\in \catlc$ and
$I$ is a negligible subset of $U$.

Let $f:X\longrightarrow Y$ be a continuous function between two
locally connected spaces.  We say $f$ is
{\em diffuse\/} if $f$ pulls back each negligible subset element of
$Y$ to a negligible subset elements of $X$.
\end{definition}

\begin{corollary}
\label{zcor1}  Let $X\in \catlc$.  Let $\cal V$ be an open cover of $X$.  Let
$(U,I)$ be a subset element of $X$.  If
$(V\cap U,V\cap I)$ is negligible for each
$V\in \cal V$, then $(U,I)$ is negligible.
\end{corollary}

\begin{corollary}
\label{ycor1}
Let $X,Y\in \catlc$, and let
$f:X\longrightarrow Y$ be a function.  Let $\cal V$ be an
open cover of $X$.  If, for each $V\in \cal V$, the restriction of $f$ to
$V$ is a diffuse function $V\longrightarrow Y$, then
$f$ is a diffuse function.
\end{corollary}

\begin{pfproof}
Trivial.
\end{pfproof}

We are ready to define the
{\em categories\/} to which we will apply the machinery of this
paper's predecessors.  Trivially,
\cattop -isomorphisms are diffuse and composition of diffuse
continuous functions are diffuse.

\begin{definition}
Let \catgen~denote the subcategory of \cattop~whose object
class is \catlc~and
whose morphism class consists of all diffuse continuous functions
between such objects.  Define the subcategory
{\bf HGen} (respectively, {\bf HGen}$^{+}$) of \catgen~to have
the class of all Hausdorff (respectively, locally Hausdorff) locally
connected spaces for objects, and the class of all
\catgen -morphisms between such objects for morphisms.
\end{definition}

The category \cattop~has more than mere objects and morphisms.
Each object supports a canonical Grothendieck topology.  It is not hard
to lift this
topology---perhaps system of topologies is a better phrase---to
\catgen , {\bf HGen} and {\bf HGen}$^{+}$.

\section{Pseudogeometric Topologies}
\label{geometric}
For this section,
\begin{onecond}
\label{ycond1}
\item	let \catc~be \catgen , {\bf HGen} or {\bf HGen}$^{+}$.
\end{onecond}
We shall define on \catc~a ``canonical'' topology, which we shall call the
{\em pseudogeometric topology\/}.  This topology has many special properties,
as abstracted in \cite{LG}, \cite{O1} and \cite{O2}.  The formulations in
these papers are non-standard.
As the reader may not be familiar with the language, we go through
verification carefully.  Each individual step is rather simple.  Difficult
points have been addressed in earlier papers, in general form.
Hopefully, work with the explicit categories of this paper will
illustrate the ethereal machinery of its predecessors.

Definition of a categorical topology, in the sense of \cite{LG},
begins with a selection of a class of special morphisms.  Let
$Sub$ denote the class of open embeddings (that is, open injections in
\cattop ) whose domain and codomain are in \catc .  Note that if
$u:U\longrightarrow X$ is an open embedding and $X\in \catc$, then
$U\in \catc$ and $u$ is a
\catc -morphism.  We claim that, in the language of \cite[(2.1)]{LG},
$Sub$ is a
{\em universe of formal subsets.}  This amounts to three conditions:
\begin{pflist}
\label{ylist1}
\item	$Sub$ contains all  \catc -isomorphisms,
\item	composition of members of $Sub$ are in $Sub$,
\item	each member of $Sub$ is a pullback base, and every pullback
of a member belongs to $Sub$.
\end{pflist}
The first two conditions are self-evident.  The third involves a subtle point.

Let $f:X\longrightarrow Y$ be a \catc -morphism and let
$u:U\longrightarrow Y$ be an open embedding.  Choose
$V$ to be the inverse image, under $f$, of the image of $u$.  Let
$v:V\longrightarrow X$ be subset injection, and let
$g:V\longrightarrow U$ be the unique function such that
$u\circ g=f\circ v$.  Assign to
$V$ the subset topology.  Then
$(V;v,g)$ is a pullback of $u$ along $f$ with respect to the category
\cattop !  More importantly, $U,V\in \catc$ and each of $u$, $v$ and $g$
is diffuse.  This suggests that the triple might also be a pullback with
respect to \catc .  Indeed,
$(V;v,g)$ is a pullback $f^{-1}u$ in \catc .  Proof relies on the following

\begin{lemma}
Assume (\ref{ycond1}).  Let $f:X\longrightarrow V$ be a
continuous function between members of {\bf LC}.  Let
$v:V\longrightarrow Y$ be an open embedding into another member
of {\bf LC}.  Then $v\circ f$ is diffuse
if and only if $f$ is diffuse.
\end{lemma}

\begin{pfproof}
Trivial.
\end{pfproof}

Categorical pullbacks along open embeddings now have an explicit
characterization.  Condition (\ref{ylist1}.c) follows immediately.

It is unusual for a pullback of a \catc -morphism, as defined in \cattop ,
to serve as a pullback in \catc .  We shall see later that \catc~is not
closed under arbitrary fibered product.  We shall struggle with morphisms
which are pullback bases but whose pullbacks in \catc~disagree with their
pullbacks in \cattop .

Having chosen a suitable $Sub$, we need a notion of cover.
We say that a non-empty cone $S$ of open embeddings into an
object $X\in \catc$ {\em covers} if
$X=\bigcup _{s\in S}~$Im$(s)$.
Formally, we choose $Cov$ to be the class of all non-empty cones of
open embeddings with this property, and must verify that $Cov$
satisfies the conditions of \cite[(2.9)]{LG}.  Most are obviously true.
\begin{pflist}
\label{ylist2}
\item	Each $S\in Cov$ is a non-empty cone.
\item	If $S\in Cov$ and $T$ is a cone of open embeddings which
	contains $S$, then $T\in Cov$.
\item	If $b$ is a \catc -isomorphism, then $\{ b\} \in Cov$.
\item	If $S\in Cov$, and if $\theta (s)\in Cov$ is a cover of
	\domk{s} for each $s\in S$, then
	$$
		\{ s\circ u~:~ s\in S, u\in \theta (s)\} \in Cov.
		$$
\item	If $S\in Cov$ is a cover of $Y\in \catc$ and if
	$f:X\longrightarrow Y$ is a \catc -morphism, then, for any choice of
	pullbacks, the set $\{ f^{-1}s~:~s\in S\}$ belongs to $Cov$.
\item	If $f\in Sub$ and $\pi _{1}$ and $\pi _{2}$ are the canonical projections
	$f\times _{cod~f}f\longrightarrow \domk{f}$, then
	$\{ \pi _{1}\} ,\{ \pi _{2}\} \in Cov$.
\end{pflist}
Condition (\ref{ylist2}.e) relies on the explicit description
of pullbacks.  Condition  (\ref{ylist2}.f) is less subtle; since each
$f\in Sub$ is monomorphic, the two projections are isomorphisms!

We refer to $(Sub,Cov)$ as the
{\em (canonical) pseudogeometric topology.}  Let us consider
terminology and properties.

In earlier works, the author refers to members of $Sub$ as
{\em formal subsets.}  A formal subset $b$ for which
$\{ b\} \in Cov$ is called a
{\em covering morphism.}  In the present, explicit, context, we
continue to refer to these key morphisms as open embeddings.  In this
topology, a morphism is a covering morphism if and only if it is a
\catc -isomorphism.

The topology meets the
{\em smallness condition\/} \cite[(2.11)]{LG}.  This is the categorical
name for the observation that, for each $X\in \catc$,
\begin{pflist}
\item	every open embedding into $X$ is $\catc /X$-isomorphic
	to a member of the set of subset injections of open subsets of $X$,
\item	given a family of open subsets of $X$, the issue of whether
	the family covers is set-theoretic.
\end{pflist}
Now suppose $\alpha$ and $\beta$ are two cones of open
embeddings into an object $X\in \catc$.  Suppose that for each
$j\in \dom{\alpha}$, there is an index $k\in \dom{\beta}$
such that $\alpha (j)$ factors through $\beta (k)$.  In addition,
suppose that $\alpha$ covers $X$.  Clearly, $\beta$ must cover $X$.
This property of the topology is called the
{\em flushness condition\/} \cite[Definition 2.19]{LG}.

Several categorical formulations rely on
{\em canopies.}  Let us first consider the canopy of a cover, as in
\cite[(2.21)]{LG}.  Let $X\in \catc$, and let $\theta$ be a
non-empty cone of open embeddings into $X$.  Without loss of
generality, we may assume that for each $j\in J=\dom{\theta}$,
$\theta (j):U(j)\longrightarrow X$ is subset injection.  Consider a
graph $A_{0}$, whose objects are indexed by $J\cup J^{2}$, in which
$$ \begin{array}{cl}
	A_{0}[j]=U(j) &	\mbox{for each }j\in J, \\
	A_{0}[j,k]=U(j)\cap U(k) &	\mbox{for }j,k\in J,
	\end{array}
	$$
and in which the only morphisms are the injections
$A_{0}[j,k]\longrightarrow A_{0}[j]$ and
$A_{0}[j,k]\longrightarrow A_{0}[k]$ for all choices
$j,k\in J$.  Consider the cone
$\alpha :A_{0}\longrightarrow X$ in which $\alpha (t)$ is subset
injection for each index $t$.  Then $A_{0}$ is a canopy of $\theta$, and
$\alpha$ is its canonical cone into $X$.

We claim that if $\theta$ is a cover, then $\alpha$ is a colimit.
Unwinding definitions, the colimit condition becomes:
\begin{onecond}
\item	Suppose that $Y\in \catc$ and $\{ f(j):U(j)\longrightarrow Y\}$
	is a family of diffuse continuous functions such that, for any
	$j,k\in J$, $f(j)$ and $f(k)$ agree on $U(j)\cap U(k)$.
	Then there is a unique diffuse continuous function
	$f:X\longrightarrow Y$ such that, for each $j\in J$,
	$f(j)$ is the restriction of $f$ to $U(j)$.
\end{onecond}
This is Corollary~\ref{ycor1}.  Note that if $\theta$ is a cover,
the canonical cone of any canopy for any pullback of $\theta$
(which is also a cover) is a colimit.  In the language of
\cite[Theorem~2.28]{LG}, the topology is
{\em intrinsic.}

Two of our categories are
{\em global structures.}  For that reason, we add some comments
on abstract canopies.

Let $A_{0}$ be an abstract canopy, in the sense of
\cite[Definition 3.4]{LG}, with respect to the pseudogeometric
topology on \catc .  Put
$J=\Lambda (A_{0})$.  For each $j\in J$, we have a
\catc -object $A_{0}[j]$; for each pair $(j,k)\in J^{2}$, we have an
object $A_{0}[j,k]$ and canonical open embeddings
$\rho _{1}:A_{0}[j,k]\longrightarrow A_{0}[j]$ and
$\rho _{2}:A_{0}[j,k]\longrightarrow A_{0}[k]$.  We are interested
in the issue of whether $A_{0}$ has an
{\em affinization.}  That is, whether there is an object $X\in \catc$
and a family of open embeddings $\alpha (j):A_{0}[j]\longrightarrow X$,
one for each $j\in J$, such that,
\begin{pflist}
\label{zlist2}
\item	$\{ \alpha (j)\} _{j\in J}$ is a cover of $X$, and
\item	for $j,k\in J$, $(A_{0}[j,k];\rho _{1},\rho _{2})$ is a
	fibered product $\alpha (j)\times _{X}\alpha (k)$.
\end{pflist}
In fact, there is an easy construction for $X$.

Let
$$
	X_{1} = \{ (x,j)~:~ j\in J,~x\in A_{0}[j]\} ,
	$$
and
$$
	R = \{ ((x,j),(y,k))~:~\exists z\in A_{0}[j,k]
	\mbox{~for which~} x=\rho _{1}(z)
	\mbox{~and~} y=\rho _{2}(z)\} .
	$$
For each $j\in J$, let $\beta (j)$ be the function $x\mapsto (x,j)$
from $A_{0}[j]\longrightarrow X_{1}$.  With respect to a unique
choice of topology on $X_{1}$, $X_{1}$ paired with the morphisms
$\beta (j)$ becomes a disjoint union of the family $\{ A_{0}[j]\} _{j\in J}$.

The first remark is that $R$ is an equivalence relation on $X_{1}$.
Reflexivity for $R$ relies on axiom \cite[(3.5.c)]{LG} that $A_{0}$
be a canopy.  The latter requires that, for each $j\in J$, there is a function
$\delta :A_{0}[j,j]\longrightarrow A_{0}[j]$ such that, for each
$x\in A_{0}[j]$, $\rho _{1}(\delta (x))=x=\rho _{2}(\delta (x))$.
Similarly, the symmetry property is a consequence of another
axiom  \cite[(3.5.d)]{LG} phrased in terms of existence of a
morphism.  The transitivity property comes from  \cite[(3.5.e)]{LG},
although here we need to know about pullbacks as well.  That
conditions states that, for
$i,j,k\in J$, there is a function
$$
	\omega : (A_{0}[i,j],\rho _{2})\times _{A_{0}[j]}(A_{0}[j,k],\rho _{1})
	~\longrightarrow A_{0}[i,k]
	$$
with good properties.  The relevance is as follows: suppose
$((x,i),(y,j))$ and $((y,j),(z,k))$ belong to $R$.  Take $r\in A_{0}[i,j]$
and $s\in A_{0}[j,k]$ for which
$$
	\rho _{1}(r)=x,   \rho _{2}(r)=y=\rho _{1}(s)
	\mbox{~~and~~}  \rho _{2}(s)=z.
	$$
Then $(r,s)$ represents a member of
$A_{0}[i,j]\times _{A_{0}[j]}A_{0}[j,k]$, and
$$
	\rho _{1}(\omega (r,s))=\rho _{1}(r)=x
	\mbox{~~and~~}   \rho _{2}(\omega (r,s))=\rho _{2}(s)=z
	\Rightarrow ~~((x,i),(z,k))\in R.
	$$
Three of our categorical formulations are no more than the axioms of an
equivalence relation!

We can now define a quotient space (with canonical projection)
$q:X_{1}\longrightarrow X(=X_{1}/R)$ and define
$\alpha (j)=q\circ \beta (j)$ for each $j\in J$.  It remains to
check (\ref{zlist2}.a,b).

Once we show that each $\alpha (j)$ is an open injection, then
(\ref{zlist2}.a) is a tautology.  First, fix $j\in J$.  Each projection
$A_{0}[j,j]\longrightarrow A_{0}[j]$ is injective, and composes with
the diagonal $\delta$ to get the identity function.  Consequently, each
projection, and $\delta$, is an isomorphism.  It follows that distinct
members of $A_{0}[j]$ are not equivalent mod$(R)$.  In other words,
$\alpha (j)$ is injective.   Next, the assumption that every projection
$A_{0}[j,k]\longrightarrow A_{0}[j]$ and
$A_{0}[j,k]\longrightarrow A_{0}[k]$ is open directly implies that, for
$U\subseteq X_{1}$ an open subset, the set of all $x\in X_{1}$ which
are equivalent to a member of $U$ is also open.  Consequently, $\alpha (j)$
is an open embedding for every $j\in J$.

We know explicitly how to take a fibered product of open embeddings
in \catc .  In particular, any construction of a fibered product in
\cattop~suffices.  Thus, we may characterize
$\alpha (j)\times _{X}\alpha (k)$ as the set
$\{ (r,s)~:~\alpha (j)(r)=\alpha (k)(s)\}$ paired with a specific topology.
That the latter must be isomorphic to $A_{0}[j,k]$ is trivial.  Condition
(\ref{zlist2}.b) follows.

There is one problem:
{\em does $X$ belong to our category?}  If so, then we may conclude that
our construction is an affinization.  If \catc~is \catgen~or
{\bf HGen}$^{+}$, it is clear that $X$ always will be an object.  At this
point, we may deduce
\begin{onecond}
\label{ycond2}
\item	The categories \catgen~and {\bf HGen}$^{+}$ are closed
	under affinization.
\end{onecond}
On the other hand, it is clear that there are choices for the canopy
$A_{0}$ in {\bf HGen} for which $X$ is {\em not\/} Hausdorff.

Next, we claim that \catc~meets the
{\em CLCS\/} criterion.  This criterion has several parts.  Let $Cvm$ be
the class of covering morphisms---in our cases, the class of all
\catc -isomorphisms.  We require that $Cvm$ be a
{\em universe of layered morphisms.}  This means, firstly, that the
analogues to  (\ref{ylist1}.a,b,c) are true with $Cvm$ in place of
$Sub$; obviously, this much is true.  In addition, we require that if
$b:B\longrightarrow A$ is a \catc -morphism and if $S$ is a cover of
$A$ such that $s^{-1}b\in Cvm$ for each $s\in S$, then $b\in Cvm$.
This implication is a trivial consequence of the fact that the topology
is intrinsic (which implies that $B$ and
$A$ are colimits of the same canopy).

Suppose $b:B\longrightarrow A$ is a \catc -morphism and $S$ is a
cover of $B$ such that, for each $s\in S$,
\begin{pflist}
\label{ylist3}
\item	$b\circ s$ is an open embedding,
\item	$b^{-1}(b\circ s)$ exists, and is an isomorphism.
\end{pflist}
The last part of the CLCS criterion demands that, under (\ref{ylist3}.a,b),
$b$ must be an open embedding.  The implication requires a short paragraph.

Let $b$ satisfy (\ref{ylist3}.a,b).  We claim that $b$ is an
open embedding.  Given (\ref{ylist3}.a), it suffices to show that
$b$ is injective.  Without loss of generality, assume that $S$
consists of injections of members of a family of open subsets
$\cal U$.  Now suppose $x,y\in B$ such that $b(x)=b(y)$.  Take
$U\in {\cal U}$  such that $x\in U$.  Condition (\ref{ylist3}.b)
translates as $U=b^{-1}(b(U))$ for each $U\in {\cal U}$ .  Thus,
$y\in U$.  But, by assumption, the restriction of $b$ to $U$ is an
open embedding, which means $x=y$.

The last axiom discussed in \cite[Definition 14.1]{LG} is that
\catc~be complete (or closed) under $Cvm$. The issue is as
follows.  Suppose $J$ is a set and $A_{0}$ and $Q_{0}$ are
two canopies of type Int$(J)$ (that is, objects indexed by
$J\cup J^{2}$, morphisms between appropriate members)
and let $q\mapsto q[\alpha ]$ be a graph transformation
$Q_{0}\longrightarrow A_{0}$.  Assume that for
$j,k\in J$, the triples $(Q_{0}[j,k];\rho _{1},q[j,k])$ and
$(Q_{0}[j,k];\rho _{2},q[j,k])$ are fibered products
$(Q_{0}[j],q[j])\times _{A_{0}[j]}(A_{0}[j,k],\rho _{1})$ and
$(Q_{0}[k],q[k])\times _{A_{0}[k]}(A_{0}[j,k],\rho _{2})$,
respectively.  We call $(Q_{0},q)$ a
{\em pullback system\/} into $A_{0}$.  If {\bf Aux} is a universe of
layered morphisms and $q[\alpha ]\in ${\bf Aux} for each
$\alpha \in J\cup J^{2}$, we call it a pullback system of
{\bf Aux}-morphisms.  Our last condition is that if $(Q_{0},q)$
is a pullback system of $Cvm$-morphisms for $A_{0}$ and if
$A_{0}$ has an affinization, then $Q_{0}$ has an affinization.
Since $Cvm$ consists of isomorphisms, the implication is vacuous.
We shall look at closure for more interesting notions of layered
morphisms shortly.

We can summarize our work so far.

\begin{theorem}
With respect to their respective pseudogeometric topologies,
the categories \catgen~and
{\bf HGen}$^{+}$ are global structures, and {\bf HGen} is a local
structure, in the sense of \cite[Definition 14.1]{LG}.
\end{theorem}

\begin{corollary}
The subcategory injection functor
${\bf HGen}\longrightarrow {\bf HGen}^{+}$ is, with respect
to the pseudogeometric topologies, a
{\em globalization.}
\end{corollary}

\begin{pfproof}
The key points in the proof of the Corollary are summarized in
\cite[Remark 14.6]{LG}.
\end{pfproof}

Our three categories are linked intimately to the topological notion
of connectedness.  The next step is to show that the topological
version implies the categorical notion of connectedness developed in
\cite[(16.a,b)]{O1} and \cite[Definition 13]{O2}.  We prove that each
of our categories is
{\em topologically componentwise.}

	Let $\emptyset$ denote the empty space.  Then
\begin{pflist}
\item	for each $B\in \catc$, there is a unique \catc -morphism
	$\emptyset \longrightarrow B$,
\item	if $b:B\longrightarrow \emptyset$ is a \catc -morphism,
	then it is an isomorphism.
\end{pflist}
That is, $\emptyset$ is an empty object in the sense of
\cite[(14.a,b)]{O1}.

Let $B,C\in \catc$, and let $A$, paired with canonical injections
$b:B\longrightarrow A$ and
$c:C\longrightarrow A$, be a disjoint union, in the topological sense.
Trivially, $A\in \catc$ is a disjoint union in the categorical sense.
That is, if $Y\in \catc$ and $f:B\longrightarrow Y$ and
$g:C\longrightarrow Y$ are diffuse continuous function, then there
is a unique diffuse continuous function $h:A\longrightarrow Y$ such that
$f=h\circ b$ and $g=h\circ c$.  The two morphisms $b$ and $c$ are
open embeddings, and the pair $\{ b,c\}$ covers $A$.  The fibered product
$b\times _{A}c$ is empty.  If $f:D\longrightarrow A$ is a
\catc -morphism, then pulling back $b$ and $c$ along $f$ determines
a disjoint union structure on $D$.

The last paragraph has several tautological implications.
\begin{pflist}
\item	\catc~is componentwise,
\item	a \catc -object is connected in the categorical sense if
	and only if it is connected in the topological sense,
\item	a \catc -morphism is complemented in the categorical
	sense if and only if it is complemented in the topological sense,
\item	the pseudogeometric topology meets the first
	and last conditions in \cite[(18)]{O2}.
\end{pflist}
Let $Comp$ denote the class of complemented morphisms.

We must show that $Comp$ is a universe of layered morphisms
under which \catc~is closed.  Conditions
(\ref{ylist1}.a,b,c) are trivial.  Now suppose
$b:B\longrightarrow A$ is a
\catc -morphism and ${\cal U}$ is an open cover of $A$ such
that, for each $U\in {\cal U}$ , the restriction of $b$ to
$b^{-1}U$ is complemented.  We leave it for the reader to check that
$b$ must be complemented.

Finally, we claim that the category is closed
with respect to complemented morphisms.  In \catgen~and
{\bf HGen}$^{+}$, every canopy has an affinization, so the
claim is tautological.  Now assume \catc ={\bf HGen}.  Suppose
$(Q_{0},q)$ is a pullback system of complemented morphisms
to a canopy $A_{0}$, and
$a\sharp :A_{0}\longrightarrow A$ is an affinization in
{\bf HGen}.  Let $q\sharp :Q_{0}\longrightarrow Q$ be an
affinization in \catgen , and let $f:Q\longrightarrow A$ be the
unique morphism such that
$f\circ q\sharp [\alpha ]=a\sharp [\alpha ]\circ q[\alpha ]$
for each $\alpha \in J\cup J^{2}$.  General nonsense implies that
$f$ is a complemented morphism
{\em with respect to\/} \catgen .  It follows that the domain of $f$ is
in {\bf HGen}, which means that $q\sharp$ is an affinization
in {\bf HGen}.

\begin{theorem}
The pseudogeometric topologies of \catgen , {\bf HGen} and
{\bf HGen}$^{+}$ are topologically componentwise.  Moreover,
in each category, every object has a cover by connected objects.
\end{theorem}

\section{The Pseudo\'etale Topology}
\label{pseudo}
Let \catd~be a topologized category which is flush,
intrinsic, closed under descent and such that every formal
\catd -subset is monomorphic.  Then the definition of superopen from
\cite[Section~5]{O3} 
makes sense in \catd .  Moreover, if \catd~is topologically
componentwise, and every \catd -object has a cover by connected
objects, then we can define the
{\em derived pseudo\'etale\/} and
{\em torsorial\/} topologies on \catd~as well.  Thus, without further
comment, it follows that \catgen~and {\bf HGen}$^{+}$ support topologies
suitable for orbifolds.

Or rather, they almost do.  Earlier papers developed the machinery
necessary to generate topologies which, in turn, generate categories
of formal quotients (and of formal quotients ``pasted'' together.
We will not add to that theory.  Instead, we raise several significant
points about that construction.  In these notes, we focus on the
pseudo\'etale topology.

\subsection{Terminology}
We have provided categorical definitions for topological words
like ``discrete'', ``open'', ``finite'', etc.  In \catgen , these phrases
need not assume their standard meanings.  Indeed, \catgen~is
created specifically as a context where a discrete morphism can
have, in the traditional sense, a small set of ramification.

Discussion of terminology is hampered by the lack of a
construction for a pullback of one \catgen -morphism along
another.  However, we will make some elementary points.
When we are using a term in the general, categorical sense,
we shall prefix it with ``c-''; otherwise, words have the usual
meaning in point-set topology.  In what follows
\begin{onecond}
\label{xhyp2}
\item	  Let \catc~be either \catgen~or {\bf HGen}$^{+}$.
\end{onecond}

\begin{proposition}
Assume (\ref{xhyp2}).  Let
$b:B\longrightarrow A$ be a \catc -morphism which overlays
absolutely.  Then $b$ is surjective, in the usual sense.
\end{proposition}

\begin{pfproof}
Unfortunately, we do not have one argument that works in general.  Instead,
we offer a line of reasoning for each choice of \catc .

First, suppose \catc ={\bf HGen}$^{+}$.  Every pullback of $b$
along a formal subset also overlays absolutely.  These pullbacks
agree with the usual sense of pullback.  Hence, without loss of generality,
we may assume that
$A$ is Hausdorff.  We proceed by contradiction.

Assume $x\in A$ is not in the image of $b$.  The subset $U=A-\{ x\}$ belongs to
\catc , and injection
$\iota :U\longrightarrow A$ is  a
\catc -morphism.  Moreover, there is a
\catc -morphism $c:B\longrightarrow U$ such that $b=\iota \circ c$.  Let
$(P;p,q)$ be a self-product
$b\times _{A}b$.  Clearly $c\circ p=c\circ q$.  Since $b$ is a
colimit of the canopy of $\{ b\}$, the latter implies that $1_{A}$
factors through $\iota$, which is absurd.

Next, suppose that \catc ={\bf Gen}.  Let $C$ denote the set
$\{ 0,1\}$ with the indiscrete topology---that is, the only
open subsets of $C$ are $\emptyset$ and $\{ 0,1\}$.  Trivially, the
point-set topological product $A\times C$ belongs to \catc .
Define two functions $f,g:A\longrightarrow A\times C$ by
$f(x)=(x,0)$ and
$$
	g(x)=\left\{ \begin{array}{ll}
		(x,0) & \mbox{if $x$ is in the image of $b$,} \\
		(x,1) & \mbox{otherwise.}
	\end{array} \right.
	$$
It is routinely verified that $f$ and $g$ are \catc -morphisms.
By inspection, $f\circ b=g\circ b$.  However, $b$ is known to
be epimorphic, and so $f=g$.  This implies that $b$ is surjective.
\end{pfproof}

\begin{corollary}
Assume (\ref{xhyp2}).  Suppose that $f:B\longrightarrow W$ is a
\catc -morphism which overlays absolutely, and that
$w:W\longrightarrow A$ is an open embedding.  Then the image of
$w$ is the same as the image of $w\circ f$.
\end{corollary}

\begin{corollary}
\label{xcol4}  Assume (\ref{xhyp2}).  Then a c-open
\catc -morphism is open in the usual sense.
\end{corollary}

\begin{corollary}
\label{wcover}  Assume (\ref{xhyp2}).  A cone of
pseudo\'etale morphisms covers if and only if the union of
the images of the members of the cone equals the
common codomain.
\end{corollary}

\begin{pfproof}
In this order, the Corollaries have obvious proof.
\end{pfproof}

\begin{corollary}
\label{xcol2}  Assume (\ref{xhyp2}).  Let
$b:B\longrightarrow A$ be a c-open \catc -morphism, and let
$c:C\longrightarrow A$ be an arbitrary \catc -morphism with
the same codomain.  Suppose $x\in C$ such that $c(x)$ is in the
image of $b$.  Then $x$ lies in the image of any choice of $c^{-1}b$.
\end{corollary}

\begin{pfproof}
Let $U$ be the image set of $b$, and let $V=c^{-1}U$.  Let
$u:U\longrightarrow A$ and
$v:V\longrightarrow C$ be subset injection.  Now $b$ is composition
of $u$ with a morphism that overlays absolutely.  Hence, the pullback
of $b$ along $u$ overlays absolutely.  Since $v$ is a pullback of $u$,
it follows that the pullback of $c^{-1}b$ along $v$, which is a pullback
of $u^{-1}b$, overlays absolutely.  Consequently, the latter is surjective.
The proposition follows now from the fact that pullback along open
embeddings agrees with the usual sense of pullback.
\end{pfproof}

\begin{lemma}
\label{xlem5}  Assume (\ref{xhyp2}).  Let
$b:B\longrightarrow A$ be a c-discrete \catc -morphism.  If
$x,y\in B$ such that
$x\neq y$ and $b(x)=b(y)$, then there is a neighborhood $U$ of
$x$ which does not contain $y$.
\end{lemma}

\begin{pfproof}
If $B\in ${\bf HGen}$^{+}$, the conclusion is true for any
two distinct points in $B$.  Assume \catc ={\bf Gen}, and
suppose that $y$ is contained in every neighborhood of $x$.  Note that
if $X$ is a subset of $B$ which contains $x$ and $y$, and if $I$ is a
closed subset of $X$ which contains $y$, then $x\in I$.

Consider the function $f:B\longrightarrow B$ which is the identity
on $B-\{ x\}$ but which maps $f(x)=y$.  It follows that if
$U\subseteq B$ is an open subset, then $U\subseteq f^{-1}U$.
Thus, $f$ is continuous.  If $T$ is a subset of $B$, the
$f^{-1}T=T-\{ x\}$ or $T$.

A criterion for connectedness is needed.  We claim that
\begin{onecond}
\label{xcond1}
\item	Let $C$ be an open subset of $B$.  Then $C$ is connected
	if and only if $C-\{ x\}$ is connected.
\end{onecond}
Let $C$ be an open subset which contains $x$.  Then $y\in C$, which
implies that $C$ lies in the closure of $C-\{ x\}$.  Implication
(\ref{xcond1}) follows.

We observe next that $f$ is diffuse.  Suppose $(U,I)$ is a
negligible subset-element of $B$.  If $x\in U$, it follows
that $f^{-1}(U,I)=(U,I)$ is negligible.  If $y\notin U$, then
$f^{-1}(U,I)=(U,I)$ is negligible.  There remains the case where
$x\notin U$ and $y\in U$.  Put $U^{\ast} =U\cup \{ x\}$ and
$I^{\ast} =I\cup \{ x\}$.  Then $f^{-1}(U,I)$ is $(U^{\ast} ,I^{\ast} )$
or $(U^{\ast} ,I)$.  Remark (\ref{xcond1}) proves negligibility in either case.

Let $Y$ be the connected component of $x$ in $B$.  Then
$f_{1}$, the restriction of $f$ to $Y$, and $g_{1}=1_{Y}$,
interpreted as diffuse maps into $B$, satisfy
$$
	b\circ f_{1}=b\circ g_{1}
	\mbox{~~and~~}   f_{1}\neq g_{1}.
	$$
Consider $D=Y-\{ x\}$, and let $d:D\longrightarrow Y$ be
subset injection.  From (\ref{xcond1}), it is easy to show that
\begin{pflist}
\item	$D\in \catgen$,
\item	$d$ is diffuse.
\end{pflist}
Then $D\neq \emptyset$ and $f_{1}\circ d=g_{1}\circ d$.  But
this contradicts the discreteness of $b$.
\end{pfproof}

\begin{proposition}
\label{xprop2}  Assume (\ref{xhyp2}).  Let
$b:B\longrightarrow A$ be an c-open, c-discrete \catc -morphism
which is finite of order $n\in \bbn$.  Then, for each $x\in A$,
$b^{-1}\{ x\}$ has at most $n$ elements.  Moreover, for each
$x\in A$ for which $b^{-1}\{ x\} \neq \emptyset$, there is a connected,
c-finite, c-open, c-discrete \catc -morphism
$c:C\longrightarrow A$ with a point $z\in C$ such that,
\begin{pflist}
\item	$c(z)=x$, and
\item	for each $y\in b^{-1}\{ x\}$, there is at least one
	morphism $f:C\longrightarrow B$ for which
	$b\circ f=c$ and $f(z)=y$.
\end{pflist}
\end{proposition}

\begin{pfproof}
Fix $x\in A$ such that $b^{-1}\{ x\} \neq \emptyset$.  Suppose $C$ is
a non-empty connected \catc -object, $c\in $Mor$_{C}(C,A)$,
$z\in C$ and $f_{1},...,f_{k}$ is a finite list
$\catc /A$-morphisms $(C,c)\longrightarrow (B,b)$ such that $c(z)=x$ and
$f_{i}\neq f_{j}$ for all pairs of distinct indices $i$ and $j$.  Suppose
$y\in b^{-1}\{ x\}$ such that $f_{j}(z)\neq x$ for every $j\in \bbn (k)$.

By the previous lemma, there is an open neighborhood $U$ of
$y$ which does not contain $f_{j}(z)$ for any $j\in \bbn (k)$.
Let $(D;d,p)$ be a pullback of the restriction of $b$ to $U$ along
$c$.  By Corollary~\ref{xcol2}, there is a point
$z^{\ast} \in d^{-1}\{ z\}$.  Let $e:E\longrightarrow D$ be
the identity map on the connected component of $z^{\ast}$ in
$D$.  Since $b$ is discrete,
$f_{i}\circ d\circ e\neq f_{j}\circ d\circ e$ for any
two distinct indices $i$ and $j$.

Let $u:U\longrightarrow B$ be subset injection.  For each index $j$,
$$
	b\circ f_{j}\circ d\circ e=
	c\circ d\circ e=b\circ u\circ p\circ e,
	$$
Put $q=u\circ p\circ e$.  By choice,
$q(z)\neq \{ f_{j}\circ d\circ e\} (z^{\ast} )$.  (We do not
claim that $q(z)$ actually equals $y$, however.)  Replacing
$c:C\longrightarrow A$ and $z$ by
$c\circ d\circ e:E\longrightarrow A$ and $z^{\ast}$,
we can now receate the hypothesis of this construction
but with an indexed list, of morphisms, of length $k+1$.

The properties of c-openness, c-finiteness and c-discreteness
are preserved by pullback and composition with complemented
morphisms.  Proof of our proposition is
now a trivial consequence of the above construction.
\end{pfproof}

\begin{corollary}
\label{xcol3}  Assume (\ref{xhyp2}).  Let
$b:B\longrightarrow A$ be a c-finite, c-discrete, c-open
\catc -morphism.  Let
$c:C\longrightarrow A$ be a \catc -morphism with the same
codomain, and let $(P;p,q)$ be a fibered product
$c\times _{A}b$.  Suppose $x\in B$ and $y\in C$ such that
$b(x)=c(y)$.  Then there is $z\in P$ such that $p(z)=y$ and $q(z)=x$.
\end{corollary}

\begin{pfproof}
Let $x$ and $y$ be as hypothesized.  Let $w=b(x)=c(y)$.  By
Proposition~\ref{xprop2} and Lemma~\ref{xlem5}, there is an
open neighborhood $U$ of $x$ which contains no other members
of $b^{-1}\{ w\}$.  In a canonical sense, $q^{-1}U$ is a
pullback of the restriction of $b$ to $U$ along $c$.  Thus, there is
$z\in q^{-1}U$ such that $p(z)=y$.  The only possible value of
$q(z)$ is $x$.
\end{pfproof}

\begin{remark}
A particular self-product may illustrate the previous Corollary.
Let $n\in \bbn$, let $G$ be the multiplicative group of complex
$n$-th roots of unity, and let $b$ be the function $b(z)=z^{n}$ from
$\bbc \longrightarrow \bbc$.  Let $X$ be the disjoint
union of $n$ copies of \bbc , indexed by $G$.  Define $\pi _{1}$
and $\pi _{2}$ on $X$ by, for each
$\omega \in G$, letting $\pi _{1}$ be the identity map and
$\pi _{2}$ be multiplication by
$\omega$ on the $\omega$-th copy of \bbc~in $X$.
In \cattop , the self-product
$b\times b$ derives from $(X;\pi _{1},\pi _{2})$ if we identify
the 0's of the copies of \bbc~as one point.  In \cattop , for each
pair $(r,s)$ of complex numbers for which
$r^{n}=s^{n}$, there is a unique member $w$ of the product such
that $r=\pi _{1}(w)$ and $s=\pi _{2}(w)$.  However, although
the underlying object of the product in \cattop~belongs to
{\bf HGen}, the two projections are not diffuse; the origin in
\bbc , a negligible set, pulls back to a non-negligible set.  With
respect to {\bf HGen}, $(X;\pi _{1},\pi _{2})$ is the self-product.
To get a product in the context of diffuse functions, we must
allow for points which
{\em cannot be distinguished\/} by $\pi _{1}$ and $\pi _{2}$.

In our example, we can say that for a pair $(r,s)$ for which
$r^{n}=s^{n}$, there is
{\em at least one\/} point in the product whose first and second
projection are $r$ and $s$, respectively.  Corollary~\ref{xcol3}
assures us that, at least for products with a c-finite,
 c-discrete, c-open morphism, any pair of points
$(x,y)$ which are sent to the same image do arise as first and
second projections of {\em something\/} in the product.
\end{remark}

\subsection{The Pseudo\'etale Topology
is {\em not\/} closed under Descent}
\label{nonquot}
Corollary~\ref{xcol3} is more important than it may appear at first.

Suppose $A_{0}$ is a canopy in \catc~(with respect to the
pseudo\'etale topology).  Put $J=\Lambda (x)$.
Let $X$ be the (topological) disjoint union of the sets
$\{ A_{0}[j]\} _{j\in J}$; for each index $j$, let $\beta [j]$ be the
canonical injection
$A_{0}[j]\longrightarrow X$.  Define a relation $R$ (signified by
$\sim$) on
$X$ by, for all appropriate choices of parameters,
$\beta [j](x)\sim \beta [k](y)$ if and only if there is
$z\in A_{0}[j,k]$ such that $\rho _{1}(z)=x$ and $\rho _{2}(z)=y$.
Reflexivity and symmetry of $R$ follow trivially.  However,
transitivity requires Corollary~\ref{xcol3}.  That is, given indices
$i,j,k\in J$, $s\in A_{0}[i,j]$ and $t\in A_{0}[j,k]$ for which
$\rho _{2}(s)=\rho _{1}(t)$, we need the existence of
$$
	w\in (A_{0}[i,j],
	\rho _{2})\times _{A_{0}[j]}(A_{0}[j,k],\rho _{1})
	$$
such that $s$ is the first projection of $w$, and $t$ is the second.

Since $R$ is an equivalence relation, we can define a quotient function
$q:X\longrightarrow Q$ for it in \cattop .  By Corollary~\ref{xcol4}, all
morphisms of the canopy are open in the traditional sense.  It
follows that $q$ is an open function.  Consequently, $Q\in \catgen$.

We claim that {\em if\/} $A_0$ has an affinization, then the specific
cone $j\mapsto q\circ \beta [j]$ is also an affinization.  To see
this, assume that $\alpha :A_{0}\longrightarrow A$ is an
affinization.  There is a unique continuous function
$f:X\longrightarrow A$ such that
$$
	\alpha [j]=f\circ \beta [j]
	\mbox{~~~for each~}j\in J.
	$$
Recall that $\alpha [j]$ is known to be open.  It follows that
$f$ is an open function.
Our conclusion amounts to two conditions:
\begin{pflist}
\label{wiso}
\item	$f$ is surjective, and
\item	for $r,s\in X$, $f(r)=f(s)$ if and only if $r\sim s$.
\end{pflist}
Condition (\ref{wiso}.a) follows from Corollary~\ref{wcover}
and the requirement that $j\mapsto \alpha [j]$ be a cover.

Let $j,k\in J$, $x\in A_{0}[j]$ and $y\in A_{0}[k]$.  We must show
that the condition
\begin{equation}
\label{weq1}
	\beta [j](x)\sim \beta [k](y)
\end{equation}
is equivalent to
\begin{equation}
\label{weq2}
	f(\beta [j](x))=f(\beta [k](y))~~\Leftrightarrow
	\alpha [j](x)=\alpha [k](y).
\end{equation}
Since $\rho _{1}\circ \alpha [j]=\rho _{2}\circ \alpha [k]$
on $A_{0}[j,k]$, statement (\ref{weq1}) implies (\ref{weq2}).
Conversely, assume (\ref{weq2}).  Another aspect of
our affinization is that $(A_{0}[j,k],\rho _{1},\rho _{2})$
must serve as $\alpha [j]\times _{A}\alpha [k]$.  In this
context, Corollary~\ref{xcol3} implies (\ref{weq1}).

It would be nice if the previous paragraphs were preparation
for proof that the pseudo\'etale topology of \catgen~is closed under
descent.  Alas, this is not the case.  Problems with quotients are
well-known; see, for example,
\cite{GS} and \cite{WT}.  What follows is our spin on some
well-known observations.

Let $V$ be a real vector space, let $\sigma$ be a non-identity
linear automorphism of $V$ such that $\sigma ^{2}=1_{V}$, and let
$G$ be the 2-group $\{ 1_{V},\sigma \}$.  Consider
$V$ modulo $G$.  Actually, we ask a more restrictive question.
We want to know whether a quotient $G\backslash V$ exists
which is pseudo\'etale.  That is, the quotient must also have
properties of finiteness and discreteness.

Let $V^{G}$ be the disjoint
union of two copies of $V$, let $\pi _{1}:V^{G}\longrightarrow V$
be the identity map on each copy and let
$\pi _{2}:V^{G}\longrightarrow V$ be the identity on one component and
$\sigma$ on the other.  Let $\Gamma$ be the graph consisting of
$V^{G}$, $V$ and $\{ \pi _{1},\pi _{2}\}$.  Existence of a good
quotient is equivalent to two requirements:
\begin{pflist}
\item	$\Gamma$ is a canopy,
\item	$\Gamma$ has an affinization.
\end{pflist}
Each point merits comment.

Let $W$ be the fixed point set of $\sigma$.  It has been observed
several times in previous papers (see \cite[Section 8]{O3}
that $\Gamma$ is a canopy provided that the equalizer
of $\{ 1,\sigma \}$ is the empty object.  The latter conditions
means that there must {\em not} be a non-empty \catgen -morphism
$d:D\longrightarrow V$ such that the image of $d$ lies in $W$.
This requires verification.

Let us suppose that a troublesome morphism
$d:D\longrightarrow V$ exists, and try to
reach a contradition.  We may assume that the domain $D$ is
connected.  If $W$ is of codimension 2 or more in $V$, then $W$
is negligible in $V$.  Consequently, $D=d^{-1}W$ is negligible in
$D$, an absurdity.  Next, suppose $W$ has codimension 1 in $V$.  Then
every affine hyperplane of $W$ is negligible in $V$ but divides
$W$ into two connected components.  Now the inverse image
of a hyperplane $H$ in $D$ must be negligible.  It follows that $d$
maps into one of the two half-spaces of $W$ bounded by $H$.
Consequently, the image of any three points in $D$ must lie on
a(n affine) line.  Hence, the entire image of $D$ lies on a line!

At this point, our argument hits a twist.  Suppose that $V$ has
dimension at least 2, so that a point in $V$ is a negligible set.
Using points to divide the image of $D$, we deduce that the image
of $D$ can not be bigger than a single point, which leads to a
contradiction.  However, if
$V$ has dimension 1 and $W=\{ 0\}$, the subset inclusion
$d:\{ 0\} \longrightarrow V$ is the unwanted equalizer!  Indeed,
in a one-dimensional real manifold, there are no non-empty
negligible subsets.  Therefore, every continuous function into
a one-dimensional real manifold is diffuse!

\begin{remark}
\label{wrem2}
The argument shows that a non-empty \catgen -morphism
into a topological manifold of dimension $n>1$ can not map
into any proper submanifold.
\end{remark}

Hereafter, assume $V$ has dimension $\geq 2$.
Consider the topological quotient
$q:V\longrightarrow G\backslash V$.
The initial remarks of this subsection tell us that if an
affinization exists, then $q$ is one.  In fact, the canopy is
so simple that it suffices to show that $q$ is pseudo\'etale.
Alas, verification that $q$ is pseudo\'etale is non-trivial.  Later,
we prove Theorem~\ref{xtheo2} which states that $q$ has the
desired properties if $W$ has codimension $\geq 2$.  However,
in the other case, it is simple to see that $q$ can not be
pseudo\'etale.

Suppose $W$ has codimension 1.  By inspection, $q(W)$
is a negligible subset of $V/G$.  Yet, $W=q^{-1}q(W)$ is not
negligible in $V$.  Hence, $q$ is not even diffuse!

\begin{remark}
Assume that $W$ has codimension 1.  We have proved that a
good quotient $G\backslash V$ does not exist in \catgen .
However, sometimes the quotient can exist in some
category.

When $V=\bbr$, we encounter an intractable problem.
The morphism $d:\{ 0\} \longrightarrow V$ has the property that
$d=\sigma \circ d$.  Enlarging the category can not change this
equation, which rules out existence of a Galois morphism
$q:V\longrightarrow Q$ with Galois group $G$.

Now assume that dim~$V\geq 2$.  Since $\Gamma$ is
a canopy, a suitable $G\backslash V$ will exist in
$\catgen ^{+}$, as discussed in the next section.  The topological
quotient is inadequate.  The problem will be studied in detail in
Section~\ref{whatis}, but we give a synopsis here.

It is possible to find \catgen -morphisms
$f,g:D\longrightarrow V$ and $x_{0}\in D$ such that
\begin{pflist}
\item	$g(x)\in \{ f(x),\sigma (f(x))\}$ for each $x\in D$,
\item	in any neighborhood of $x_{0}$, there exist $y,z$ such that
	$$\begin{array}{ll}
		g(z)=f(z)\neq \sigma (f(z)) &\mbox{ and } \\
		g(y)=\sigma (f(y))\neq f(y) &
		\end{array}$$
\end{pflist}
Let $q:V\longrightarrow G\backslash V$ be a topological quotient.
In the topological sense, $q\circ f=q\circ g$.  If $q$ were an
affinization, then there would be a product map
$\delta :D\longrightarrow V^{G}=q\times q$ such that
$f=\pi _{1}\circ \delta$ and $g=\pi _{2}\circ \delta$.  Let
$E$ be the connected component of $x_{0}$.  Then, for all
$x\in E$, we would have one of the identities
$$
	f(x)=g(x) \mbox{~~or~~}f(x)=\sigma (g(x))
	$$
true for {\em all} $x\in E$.  Yet, this is not the case.

We have assumed too much structure, and reached a contradiction.
The reason
is that, for $q\sharp :V\longrightarrow Q$ the actual quotient
in an enlarged category,
$q\sharp \circ f$ and $q\sharp \circ g$ are not the same.
\end{remark}

\subsection{Expansions}
Suppose \catc~is a category of orbifolds of active interest.
Attaching a topological space to each object amounts to
introducing a continuous functor
$\Gamma :\catc \longrightarrow ${\bf Gen}.  Extending the
topological model to an enlargement of
\catc , via the plus functor, means lifting $\Gamma$ to a
new continuous functor $\Gamma ^{+}$ on \catcpl .  Unfortunately,
since \catgen~is not closed under descent, there may not be an
extension to $\catcpl \longrightarrow \catgen$.

It is necessary to have $\Gamma ^{+}$ defined on
all members of \catcpl , which means we need a codomain.
The ``obvious'' choice is $\catgen ^{+}$.  Indeed, if $\catgen ^{+}$
exists, the extension theorem follows from work in this series.
The problem is that we can not yet prove that the plus
construction applies to \catgen !

The problem is formal, not substantial.  As observed in
\cite[Section~3]{O3}, 
the only obstruction is to show that the class of global
classes between two $\catgen ^{p}$-objects is representable
by a set.  In fact, we make the

\begin{conjecture}
The pseudo\'etale topology for \catgen~meets the smallness
axiom for a categorical topology.  The torsorial topology
does not meet the smallness condition, but does satisfy
\cite[(11)]{O3}. 
\end{conjecture}

\noindent
After all, every pseudo\'etale morphisms into $X\in \catc$ is
in the set of finite-to-one maps to $X$.  Unfortunately, the
conditions used to define the pseudo\'etale topology are not
framed in set-theory.

It is not difficult to get around the immediate problem.  Let
$S$ be an infinite set, and let $\catgen [S]$ be the subcategory
of objects which admit a cover by open sets of cardinality less than
or equal to the cardinality of $S$.  Assign to $\catgen [S]$ the
obvious restriction of the topology.  All of the arguments in
\catgen~in this paper apply equally well $\catgen [S]$.  The
difference is that $\catgen [S]$ is the globalization of a
{\em small\/}
local structure---specifically, of the subcategory of
\catgen~objects entirely of cardinality less than or equal
to the cardinality of $S$!  In practical examples, $S=\bbr$
generates a category containing the desired topological models.

Technically, the derived pseudo\'etale topology of
\catgen~restricted to $\catgen [S]$ need not be the derived
pseudo\'etale topology of $\catgen [S]$.  That is, at its face,
it is possible that a $\catgen [S]$-morphism $f$ might be
pseudo\'etale in $\catgen [S]$ but not in \catgen .  This seems
unlikely.  However, until we have better control over pullbacks,
we can only

\begin{conjecture}
Let $S$ be an infinite set, and let $f$ be a $\catgen [S]$-morphism.
\begin{pflist}
\item	$f$ is torsorial in the category $\catgen [S]$ if and only if it
	is torsorial in \catgen .
\item	$f$ is pseudo\'etale in $\catgen [S]$ if and only if it
	is pseudo\'etale in \catgen .
\end{pflist}
\end{conjecture}

\noindent
Even without the conjecture, we have a context which allows for
arbitrary repition of the plus construction.

\subsection{Functors and Pseudo\'etale Topologies}
\label{contrast}
Let $S$ be a set, let \catm~be a pseudogeometric topologized
category and let
$\Gamma :\catm \longrightarrow \catgen [S]$ be a continuous
functor with respect to the pseudogeometric topology on
$\catgen [S]$.  For example, $\Gamma$ might be the forgetful functor
on a category of manifolds \catm , which might be differentiable,
Riemmanian, analytic, complex, etc.  It is not necessarily true that
$\Gamma$ is continuous with respect to the derived
pseudo\'etale topologies.

In practice, it seems reasonable.  However,  the author
suspects that proof of continuity in a specific case depends
on the particulars of \catm .  An old example illustrates
the problem.  Let $n\in \bbn$, $n>1$, and let
$g:\bbc \longrightarrow \bbc$ be the function $g(z)=z^{n}$.
We shall prove in Section~\ref{fullquot} that $g$ is
pseudo\'etale in the context of \catgen .  In the sense
of real manifolds, $g$ does not even overlay absolutely.
Yet, in the category of complex analytic manifolds with
diffuse morphisms, $g$ appears again to be pseudo\'etale.  What
we expect to be true is that if $h$ is pseudo\'etale in the sense
of real manifolds, or any other kind of structure, it must
pseudo\'etale in \catgen~{\em and\/} if $h$ is pulled back
along another morphism $f$,
then the underlying space of the pullback is the pullback
$f^{-1}h$ in \catgen .

It may seem odd that the author introduce a notion whose
compatibility with continuous functors is questionable.  There
is logic behind the formulation of pseudo\'etale morphisms.
Our thesis has been that the parameters of formal objects
are set at he categorical level.  Even the study of $g$ in
the above paragraph suggests that {\em context\/} rather
than intrinsic nature determine whether singularities
prevent good structure.

For all its formalism, the definition of being pseudo\'etale
is negative.  It is motivated by a question: Given a morphism
$b$, how can one determine that, even in an enlarged
category, $b$ will not obey the rules of manipulation that
we need?  Rather than turn the question on concrete issues,
such as the shape of singular sets, the author started with a
list of categorical manipulations to be allowed and tried
to define a pseudo\'etale morphism as anything for which
those manipulations did not lead to contradiction.

The advantage of such a definition is that a category like
\catgen~comes with an intrinsically defined notion of
pseudo\'etale morphism, even before singular sets are
studied.  Of course, in order to show that a morphism
with singularities actually meets the abstraction requires
work.  Only after Sections~\ref{prod} and \ref{fullquot}
will we have interesting examples.

It is not hard to modify the definition in the context of a model.
Let $\Gamma :\catm \longrightarrow \catgen [S]$ be as above.
As usual, interpret it as assigning an underlying base space to each
\catm -object.  If we want to enlarge \catm~in a manner that
assures that $\Gamma$ lifts, a modified topology can be used.
Define a \catm~morphism $b$ to be $\Gamma$-pseudo\'etale if and
only if
\begin{pflist}
\item	$b$ is pseudo\'etale in \catm ,
\item	$\Gamma (b)$ is pseudo\'etale in $\catgen [S]$, and
\item	$\Gamma$ preserves all pullbacks along $b$.
\end{pflist}
The $\Gamma$ will be continuous with respect to the
corresponding topology on \catm~and the pseudo\'etale topology
of $\catgen [S]$.

The author's guess is that, in examples of interest, continuity
will be true but very hard to prove.  A situation in which a morphism
can be pseudo\'etale but not $\Gamma$-pseudo\'etale would be
exceptionally interesting.

One key reason for pseudo\'etale topologies is the need for
a topology of finite-to-one morphisms to enable the plus
construction to produce quotients.  There is a very easy way to
define such topologies.  Let \catc~be a pseudogeometric
topology (one in which formal subsets are monomorphic).  A
\catc -morphism $b$ is called a
{\em local subset} if there is a cover
$S$ of \domk{b} such that $b\circ s$ is a formal subset for
every $s\in S$.  Let $FDL$ be the class of all local subsets which
are finite and discrete.  With an obvious notion of cover, $FDL$ becomes
a topology which meets the axioms of a pseudo\'etale topology.
We call it the {\em elementary finite-to-one topology}.
It is not as rich as the  {\em derived\/} pseudo\'etale topology.

For example, in the study of $G$ and $V$ from
Subsection~\ref{nonquot}, we would conclude that
$G\backslash V$ {\em never\/} exists in \catgen , but does exist
in $\catgen ^{+}$ whenever dim $V\geq 2$.  Because of the
simplicity of the topology, the plus functor introduces so many
new objects that even when $W$ has codimension $\geq 2$,
the topological quotient $q:V\longrightarrow G\backslash V$
loses the quotient property in the enlarged category.

In the real case, it is the elementary finite-to-one topology
that is in common usage.  For that reason, we suspect the
following to be true

\begin{conjecture}
Let \catm~be the subcategory of $C^{n}$-manifolds (where
$n\in \bbn$ or $n=\infty$) consisting of all manifolds but only
diffuse $C^{n}$-functions.  Then the pseudo\'etale and elementary
finite-to-one topologies agree.
\end{conjecture}

\subsection{Missing Things}
There are two kinds of objects that our theory does not allow
for.  One is not of great technical concern, the other is more
serious.

In limiting our categories' version of morphism, we have effectively
abolished any notion of sub-object.  Closed subsets can not be
translated into morphisms in our context.  Products are rare,
which means that an abstract form of the Rank Theorem
can not be formulated.

More serious is the lack of group objects.  Aside from
(categorical analogues to) finite
groups), there are virtually
no abelian group object.  This impacts on our ability to define
cohomology for our very formal objects.

Each object supports
a topology in the sense of Grothendieck.  Consequently,
cohomological theories arise naturally from any ``sheaf'' into
the category of abelian groups.  In our language,
{\em functors of sections\/} are used in place of sheaves,
and our theory includes
existence and uniqueness of liftings for such functors.

In many theories, for $A$ an abelian group object, the
functor Mor$(*,A)$ is a functor of sections and begats
cohomology.  Sometimes, one can adapt traditional group objects to the
task.  Suppose $A$ is an abelian topological group.  Although
$A$ may not exist in \catgen , we can define the functor
$\Phi =$Mor$_{Top}(*,A)$.  Because $A$ is outside \catgen ,
it does not follow immediately that $\Phi$ is a functor of
sections.  If, however, it is, then it will lift to all expansions
of \catgen~via the plus construction.  Verification that
$\Phi$ is a functor of sections can be simplified by using
\cite[Proposition 10]{O3} 

\subsection{What is the Hausdorff Condition?}
We have tracked Hausdorff objects because of precedent.  Historically,
objects of interest have some sort of separation property.  Although
the Zariski topology of a scheme is not Hausdorff, there is a
weaker notion of ``separated''.  Unfortunately, the author is
unaware of a separation property that can be characterized at
the universal level.

It is well-known that a descent of separated (or Hausdorff) things
can produce an unseparated object.  However, it seems that global
objects are of less interest unless they are separated (and, typically,
every local object is separated.)

The author has ideas on this topic.  Indeed, we have a notion of
separation, framed in the context of a topologically componentwise
category, with respect to which the following is true:
\begin{onecond}
\item	Let $b:B\longrightarrow A$ be a morphism and let $S$ be
	a cover of $B$ such that $b\circ s$ is discrete for each
	$s\in S$.  If $B$ is ``separated'', then $b$ is discrete.
\end{onecond}
However, before a definition can be made profitably,
a more precise understanding of what separation should entail
is needed.

We have no further results specific to {\bf HGen} or
{\bf HGen}$^+$.

\section{Representability of a Functor}
\label{prod}
Let $X\in \cattop$ and let \cate~be a family of closed subset
elements of $X$.  Define a contravariant functor $F$, denoted by
{\em Neg(X,\cate ),\/} from \catgen~to the category of sets as
follows.  For $B\in \catgen$, let
$F(B)$ be the set of continuous functions $f:B\longrightarrow X$
such that $f^{-1}E$ is negligible for every $E\in \cate$.  For
$b:B\longrightarrow C$ a \catgen -morphism, define
$F(b):F(C)\longrightarrow F(B)$ by $f\mapsto f\circ b$.

Various issues reduce to representability of a functor this kind.
For example, let $A\in \catgen$ and let
$\{ b_{j}:B_{j}\longrightarrow A\}_{j\in J}$ be a
list of members of $\catgen /A$.  Let
$(P;\{ \pi _{j}\}_{j\in J})$ be a fibered
product of the family in the topological sense (that is, in
$\cattop /A$).  Let
\cate~consist of all subset elements in $P$ which are a pullback of a
negligible subset element by one of the projections $\pi _{j}$.  In an
obvious sense, a representative for Neg$(P,\cate )$ will be a product
in the category $\catgen /A$.

\begin{remark}
The above interpretation enables us to find choices for
which the functor can not be represented.  Specifically, certain products
do not exist in \catgen .  Let $X$ be a real manifold, and consider the 
usual topological $X\times X$.  It is easily
verified that both projections $X\times X\longrightarrow X$ are diffuse.
Suppose a product $P$ for $X$ with itself in \catgen~exists.  Then
there is a canonical continuous function
$P\longrightarrow X\times X$ and a canonical \catgen -morphism
$X\times X\longrightarrow P$.  By general nonsense, these functions
are continuous bijections.  Thus, $X\times X$ is a product in \catgen .
Now there must be a diagonal \catgen -morphism
$\delta :X\longrightarrow X\times X$ is whose composition with
either projection is the identity.  Tautologically, this must be the
usual diagonal embedding.  Alas, by Remark~\ref{wrem2},
that function is not diffuse.

It is easy to build examples from fibered products instead of a
self-product.  In general, if $\Gamma$ is a graph of \catgen -objects and
morphisms, and if
$\alpha :A\longrightarrow \Gamma$ is an inverse limit for $\Gamma$
such that every morphism of $\alpha$ is diffuse, then an inverse limit in
\catgen~exists only if $\alpha$ is one.  Unfortunately, diagonal maps
imply that the candidate $\alpha$ will often fail.  Indeed, it is difficult
to see how a function can have
even a self-product unless it is ``generically'' discrete in some sense.
\end{remark}

We can not advance without some non-trivial pullback bases.  In this
section, we develop one gimmick, which will ultimately allow us to
prove something about group quotients.  From a pair $(X,\cate _{0})$,
we shall construct a pair
$(\Lambda ,\cate _{1})$ and a continuous functor
$\lambda :\Lambda \longrightarrow X$ such that
$f\mapsto \lambda \circ f$ induces a natural isomorphism of functors
Neg$(\Lambda ,\cate _{1})\longrightarrow $Neg$(X,\cate _{0})$.

The construction requires several pages.  Fix $X\in \cattop$.  In what
follows, we introduce some temporary terminology, for use just in the
construction.  In the present context, for each $x\in X$, let $X(x)$ denote
the set of open neighborhoods of $x$ in $X$.

Let $CSE(X)$ be the set of all closed subset elements of $X$.  In the
present context, call a subset $\catd \subseteq CSE(X)$
{\em closed\/} if the following conditions are true:
\begin{pflist}
\label{zlist3}
\item	For $U$ an open subset of $X$, $(U,\emptyset )\in \catd$.
\item	If $(U,I)\in \catd$ and
	$V$ is an open subset of $U$, then $(V,V\cap I)\in \catd$.
\item	If $(U,I)\in \catd$ and $J$ is a subset of $U-I$ for which
	$(U-I,J)\in \catd$, then $(U,I\cup J)\in \catd$.
\item	If $(U,I)\in CSE(X)$ and, for each $x\in I$ there exists an open subset
	$V$ of $U$ for which $(V,V\cap I)\in \catd$, then $(U,I)\in \catd$.
\end{pflist}
Obviously, an arbitrary intersection of closed subfamilies is closed.  Thus,
given a subset $\cate \subseteq CSE(X)$, there is a smallest closed
subfamily which contains \cate .  Refer to this as the
{\em closure of\/}  \cate , and denote it by $\cate ^{\ast}$.

Let $\cate _{0}\subseteq CSE(X)$, and let $f:Y\longrightarrow X$ be
a continuous function on a member of \catgen~which pulls back each
member of $\cate _{0}$ to a negligible subset element.  By earlier
lemmas, the family of all members of $CSE(X)$ whose pullback under
$f$ is negligible is closed.  In particular, $f$ pulls back each member of
${\cal C}_{0}^{\ast} $ to a negligible subset.

In what follows, assume a family $\cate _{0}\subseteq CSE(X)$ has
been specified, and put $\cate ={\cal E}_{0}^{\ast} $.  Then
\begin{onecond}
\item	\cate~has properties (\ref{zlist3}).
\end{onecond}
In addition, partially order \cate~be the relation that $(U,I)\leq (V,J)$
if and only if $U\subseteq V$ and $U-I\subseteq V-J$.

Fix $x\in X$.  Let $\cate (x)$ denote the subset of $(U,I)\in \cate$
such that $x\in U$.  Let
$\Lambda (x)$ denote the set of function $f$ on $\cate (x)$ such that
\begin{pflist}
\item	for each $(U,I)\in \cate (x)$, $f(U,I)$ is a connected component
	of $U-I$ whose closure contains $x$,
\item	for $(U,I),(V,J)\in \cate (x)$, $f(U,I)\subseteq f(V,J)$
	if $(U,I)\leq (V,J)$.
\end{pflist}
It is possible that $\Lambda (x)$ is empty.

Put
$$
	\Lambda =\{ (x,f):~x\in X,~f\in \Lambda (x)\} .
	$$
Define $\lambda :\Lambda \longrightarrow X$ by
$\lambda (x,f)=x$.  Now we need a topology for $\Lambda$.

Let $U$ be an open subset of $X$.  Define
$$
	\cate [U]=\{ I\subseteq U~:~(U,I)\in \cate \} .
	$$
For $(x,f)\in \Lambda$ such that $x\in U$, define
\begin{equation}
\begin{array}{l}
	N(U,x,f)=N(U,(x,f)) \\
	~~=\{ (y,g)~:~y\in U
	\mbox{~and~} f(U,I)=g(U,I)
	\mbox{~for all~} I\in \cate [U]\} .
	\end{array}
\end{equation}
We claim next that the family
\begin{equation}
\label{yeq1}
	\{ N(U,x,f)~:~x\in X,~U\in X(x),~f\in \Lambda (x)\}
\end{equation}
is a basis for a topology on $\Lambda$.

Let $(x,f)\in \Lambda$ and $U\in X(x)$.  Tautologically,
\begin{pflist}
\label{ylist5}
\item	$(x,f)\in N(U,x,f)$,
\item	$\lambda (N(U,x,f))\subseteq U$,
\item	for $(y,g)\in N(U,x,f)$, $N(U,y,g)=N(U,x,f)$.
\end{pflist}
In addition, if $V\in X(x)$, it is elementary to check that
$$
	N(U\cap V,x,f)\subseteq N(U,x,f)\cap N(V,x,f),
	$$
by the nature of connected components.  A first consequence
of these observations is that (\ref{yeq1}) is, indeed the basis
of a topology on $\Lambda$.  Hereafter, assign to $\Lambda$
this topology.  With that definition in hand, we may also deduce that
\begin{pflist}
\item	$\lambda$ is continuous,
\item	for each $x\in X$,
	$\{ N(U,x,f)~:~U\in X(x),~f\in \Lambda (x)\}$
	is a basis of open neighborhoods at $x$.
\end{pflist}
Let
$$
	\cate _{1}=\{ \lambda ^{-1}\alpha ~:~\alpha \in \cate \} .
	$$
The next objective is verification that composition with the function 
$\lambda$ determines a natural
isomorphism Neg$(\Lambda ,\cate _{1})\longrightarrow $Neg$(X,\cate )$.

To begin with, let us expand on (\ref{ylist5}.b).  Let
$(x,f)\in \Lambda$, $U\in X(x)$ and $I\in \cate [U]$.
Suppose $(z,h)\in N(U,x,f)$.  Now $f(U,I)=h(U,I)$.  By
definition, $z$ lies in the closure of $h(U,I)$.  Thus,
\begin{equation}
\label{wclose}
	z\in \overline{f(U,I)}~~~ \Rightarrow ~~~
	z\in I \mbox{~or~} z\in f(U,I),
	\end{equation}
because $f(U,I)$ is closed in $U-I$.  In other words, the image of
$N(U,x,f)$ under $\lambda$ lies in
$\overline{f(U,I)}\subseteq f(U,I)\cup I$ for each $I\in \cate [U]$.

Now suppose $Y\in \catgen$ and
$\alpha ,\beta :Y\longrightarrow \Lambda$ are continuous
function such that $\lambda \circ \alpha$ and
$\lambda \circ \beta$ pull back each member of \cate~to
negligible elements.  Suppose $\alpha \neq \beta$.  We claim
\begin{equation}
\label{xeq2}
	 \lambda \circ \alpha \neq \lambda \circ \beta .
\end{equation}
Fix $y\in Y$ for which $\alpha (y)\neq \beta (y)$.

If $\lambda (\alpha (y))\neq \lambda (\beta (y))$,
we are done.  Assume $x\in X$ such that
$$
	\alpha (y)=(x,f)   \mbox{~~and~~}   \beta (y)=(x,g)
	$$
where $f\neq g$.  Take $(U,I)\in \cate (x)$ such that
$f(U,I)\neq g(U,I)$.  Let $W$ be the connected component of $y$ in
\begin{equation}
\label{yeq2}
	\alpha ^{-1}N(U,x,f)~\cap ~\beta ^{-1}N(U,x,g).
\end{equation}
The set
$$
	J=(\{ \lambda \circ \alpha \} ^{-1}I \cup
	 \{ \lambda \circ \beta \} ^{-1}I)\cap W
	$$
is negligible in $W$.  The set $W-J$ is non-empty and
connected.  Let $w\in W$.  By choice, and by (\ref{wclose}),
$\lambda (\alpha (w))\in f(U,I)$ and
$\lambda (\beta (w))\in g(U,I)$.  But $f(U,I)$ and $g(U,I)$
are distinct connected components.  Thus,
$\lambda (\alpha (w))\neq \lambda (\beta (w))$.  We have (\ref{xeq2}).

Next, suppose that $Y\in \catgen$ and
$\gamma :Y\longrightarrow X$ is a continuous function which
pulls back members of \cate~to negligible elements.  For $y\in Y$, define
$f^{y}$ on $\cate (x)$, for $x=\gamma (y)$, as follows.  Suppose
$(U,I)\in \cate (x)$.  Let $W$ be the connected component of
$y$ in $\gamma ^{-1}U$, and let
$J=\gamma ^{-1}I$.  Then $W-J$ is a non-empty, connected open set.
Let $f^{y}(U,I)$ be the connected component of $U-I$ which contains
$\gamma (W-J)$.  Now
$$
	\overline{W}\subseteq \gamma ^{-1}(\overline{f^{y}(U,I)})
	~~\Rightarrow    	~~\gamma (y)\in \overline{f^{y}(U,I)}.
	$$
Trivially, the function $f^{y}$ belongs to $\Lambda (x)$.

Define $F:Y\longrightarrow \Lambda$ by
$F(y)=(\gamma (y),f^{y})$.  Then $\lambda \circ F=\gamma$.
Now suppose $y\in Y$, $(x,f)\in \Lambda$, $U\in X(x)$ and
$F(y)\in N(U,x,f)$.  By (\ref{ylist5}.c),
$N(U,x,f)=N(U,\gamma (y),f^{y})$.  Let $W$ be the
connected component of $y$ in $\gamma ^{-1}U$.  Now suppose
$w\in W$ and $I\in \cate [U]$.  Let
$J=\gamma ^{-1}I$, and then $W-J$ is a non-empty,
connected subset.  By definition, $f^{y}(U,I)$ and $f^{w}(U,I)$
are exactly the same!  Thus,
$F(w)\in N(U,x,f)$.  It follows that $F$ is continuous.

\begin{proposition}
\label{xprop1} Let $X\in \cattop$ and let $\cate _{0}$ be a
family of closed subset elements of $X$.  We denote the triple
$(\Lambda ,\cate _{1},\lambda )$ constructed above by
$(X,\cate _{0})\sharp$.  Now $\lambda$ induces a natural
transformation from
Neg$(X,\cate _{0})\longrightarrow $Neg$(\Lambda ,\cate _{1})$,
and this is an isomorphism of functors.
\end{proposition}

We now study $\Lambda$ in certain cases.

\begin{lemma}
\label{ylem1}  Let $Y\in \catgen$, let $(U,I)$ be a negligible
subset element of $Y$, and let $y\in U$.  Then there is a unique connected
component of $U-I$ whose closure contains $y$.
\end{lemma}

\begin{pfproof}
Let $W$ be the connected component of $y$ in $U$.  Then
$y$ is in the closure of the connected, open and closed subset
$W-I$ in $U-I$.  Since $I$ is negligible, $W-I$ is connected.
\end{pfproof}

\begin{corollary}
\label{ycor2} In the context of Proposition~\ref{xprop1}, suppose
that $V\subseteq X$ is an open subset such that,
\begin{pflist}
\label{ylist6}
\item	in the subset topology, $V\in \catlc$,
\item	for each $(U,I)\in \cate _{0}$, $(V\cap U,V\cap I)$
	is a negligible subset element.
\end{pflist}
Then the restriction of $\lambda$ to
$\lambda ^{-1}V$ is a bijection onto
$V$ whose inverse is continuous.
\end{corollary}

\begin{pfproof}
Clearly, the set of all elements $(U,I)$ such that
$(V\cap U,V\cap I)$ is negligible is closed in $CSE(X)$.
Specifically, every member of $\cate =\cate _{0}^{\ast}$
has this property.  Let
$\iota :V\longrightarrow X$ be the subset injection function.
Then there is a unique continuous function
$F:V\longrightarrow \Lambda$ such that
$\lambda \circ F=\iota$.  To finish, it suffices to show that
$\lambda$ is injective on $\lambda ^{-1}V$.  That is, given
$x\in V$, there is a unique
$f\in \cate (x)$.  But this is immediate by Lemma~\ref{ylem1}.
\end{pfproof}

\begin{lemma}
In the context of Proposition~\ref{xprop1}, if $X$ is Hausdorff,
then $\Lambda$ is Hausdorff.
\end{lemma}

\begin{pfproof}
Let $p$ and $q$ be distinct points in $\Lambda$.  We must
show that $p$ and $q$ can be separated by open neighborhoods.
If $\lambda (p)\neq \lambda (q)$, this is trivial.  Assume
$x\in X$, $p=(x,f)$ and $q=(x,g)$ for $f\neq g$.

Assume that $(U,I)\in \cate$ such that $f(U,I)\neq g(U,I)$.  Suppose
$(z,h)\in N(U,x,f)\cap N(U,x,g)$.  Then
$$
	f(U,I)=h(U,I)=g(U,I)
	$$
which would contradict choice of $(U,I)$.
\end{pfproof}

\begin{proposition}
In the context of Proposition~\ref{xprop1}, assume that
$K$ is a closed subset of $X$ and that
$V=X-K$ satisfies (\ref{ylist6}.a,b).  Assume also that
$(X,K)\in \cate _{0}$.  Then $\Lambda$ is locally connected and
$\lambda ^{-1}K$ is a negligible subset.  Moreover, if \cate~contains
all negligible subset elements $(U,I)$ in which $U\subseteq V$,
then, with respect to $\lambda$, $\Lambda$ represents
Neg$(X,\cate _{0})$.
\end{proposition}

\begin{pfproof}
For each $x\in X-K$, let $x^{\ast}$ denote the unique member of
$\lambda ^{-1}\{ x\}$.

Let $(x,f)\in \Lambda$ and $U\in X(x)$.  Put
$W=f(U,U\cap K)$.  For each $I\in \cate [U]$, $W-I$ is
connected.  Since $(U,(U\cap K)\cup I)\leq (U,U\cap K)$, we get that
$$
	f(U,(U\cap K)\cup I) = W-I.
	$$
The same reasoning applies to $(U,(U\cap K)\cup I)\leq (U,I)$.
The set $f(U,I)$ is uniquely characterized as the only connected
component of $U-I$ which intersects $W-I$.
Now let $w\in W$.  Working directly from the definition, we get that
$w^{\ast }\in N(U,x,f)$.  It follows easily that
$$
	\{ w^{\ast }~:~w\in f(U,U\cap K)\}  =
	N(U,x,f)~\cap  \lambda ^{-1}(X-K)\} .
	$$
Recall that $f(U,U\cap K)$ is non-empty.

We draw several conclusions.  First, $\lambda ^{-1}(X-K)$ must be
dense in $\Lambda$.  Second, suppose $(x,f)\in \Lambda$ and
$U\in X(x)$.  Then
$N(U,x,f)\cap \lambda ^{-1}(X-K)$ is a dense, connected subset of
$N(U,x,f)$.  Thus, $N(U,x,f)$ is connected.

We have proved that $\Lambda$ is locally connected and that
$\lambda ^{-1}K$ is nowhere dense.  In addition, for each
$r\in \lambda ^{-1}K$, we have found that, for each $U\in X(\lambda (r))$,
$N(U,r)$ is an open, connected neighborhood such that
$N(U,r)-\lambda ^{-1}K$ is also connected.  By Lemma~\ref{ylem2},
the set $\lambda ^{-1}K$ is negligible.  Consequently, $\cate _{1}$
contains only negligible subset elements.

Finally, suppose that \cate~contains all negligible subset elements
$(U,I)$ in which $U\subseteq X-K$.  If $(U,I)$ is a negligible subset
element of $\Lambda$, it is elementary to check that $(U,I)$ belongs to
$\cate _{1}^{\ast}$.  It follows that $\Lambda$, as a
\catgen -object, represents Neg$(X,\cate _{0})$.
\end{pfproof}

\begin{corollary}
\label{xcol1} Let
$b:B\longrightarrow A$ be a \catgen -morphism.  Let
$I\subseteq B$ be a negligible subset, and assume that the restriction of
$b$ to $B-I$ is, in the usual sense, an open local homeomorphism.  Then
$b$ is a pullback base in \catgen .

In fact, we can give a construction for pullbacks.  Let
$c:C\longrightarrow A$ be a \catgen -morphism.  Let
$(X;\pi _{C},\pi _{B})$ be a pullback $c^{-1}(B,b)$ with respect
to the category \cattop .  Let
$K=\pi _{B}^{-1}I$.  Let \catd~be the set of all negligible subset
elements of $X-K$, and let $\cate _{0}=\catd \cup \{ (X,K)\}$.  Let
$(\Lambda ,\cate _{1},\lambda )=(X,\cate _{0})\sharp$.  Then
$(\Lambda ;\pi _{C}\circ \lambda ,\pi _{B}\circ \lambda )$
is a pullback $c^{-1}(B,b)$, with respect to \catgen .
\end{corollary}

\begin{pfproof}
Let us justify the construction of $c^{-1}(B,b)$.  First, observe
that the restriction of $\pi _{C}$ to $X-K$ is an open local
homeomorphism.  Thus, $X-K\in \catlc$.

Let $\cal N$ be the class of negligible subset elements of $C$.
Since negligibility is a local property, it follows that for
$(U,I)\in {\cal N}$,
$((\pi _{C}^{-1}U)-K,I-K)$ is negligible.  In other words,
every subset element of
$X-K$ which belongs to $(\pi _{C}^{-1}{\cal N})^{\ast }$ is
negligible.  Just as importantly, the fact that $\pi _{C}$ is a
local homeomorphism implies that
$(\pi _{C}^{-1}{\cal N})^{\ast }$ contains
{\em every\/} negligible subset element of $X-K$.

We can find an open cover for $B-I$ such that, on each member,
$b$ restricts to an open embedding.  Consequently, there is an
open cover of $X-K$ such that $\pi _{B}$, restricted to each,
identifies with a restriction of $c$ to a subset of $C$.  By
Corollary~\ref{ycor1} the restriction of $\pi _{B}$ to $X-K$ is diffuse.

Verification that $\Lambda$ has the fibered product property
is now trivial.
\end{pfproof}

\section{Quotients and the Pseudo\'etale Property}
\label{fullquot}
We are ready to use the material of
\cite[Section~8]{O3}.  
We assume the notational conventions and theorems of that work,
with one caveat.  In that paper, the words ``open'', ``discrete'',
``finite'' and several others have categorical definitions.  In the
topological context, these terms have more standard, concrete
meaning.  In general, when using a term, we mean it in the
topological sense.  To indicate when an ambiguous term is
intended in the universal sense, we prefix it with ``c-''.

Let {\bf Aux} be the class of tuples $(G,\rho ,B,b,A)$ in which
\begin{pflist}
\item	$B\in \cattop$,
\item	$(G,\rho )$ is a finite group action on $B$, with
	respect to \cattop ,
\item	$b:B\longrightarrow A$ is an open, continuous surjection, and
\item	for $x,y\in B$, $b(x)=b(y)$ if and only if $x$ is in the
	$G$-orbit of $y$.
\end{pflist}
Note that given a group action $(G,\rho )$ on a topological space
$B$, then the standard topological quotient construction produces
such a function $b$.  In the present context, for $g\in G$ and
$x\in B$, we write $g\cdot x$ for $\rho (g)(x)$.

Let $(G,\rho ,B,b,A)\in ${\bf Aux}.  Let $U$ be the set of
$x\in B$ which have a neighborhood
$V$ such that $g\cdot V\cap V=\emptyset$ for
every $g\in G-\{ e\}$.
Obviously, $U$ is open and $G$-invariant.  In the present Section,
refer to the complement of $U$ as the
{\em upper ramification set,\/} and refer to the image of
the complement under $b$ as the
{\em lower ramification set.}  Clearly, both ramification sets
are closed in their respective spaces.  Note that the restriction of
$b$ to $U$ is a local homeomorphism.

The remainder of the section is dedicated to proof of the following

\begin{theorem}
\label{xtheo2}  Let $(G,\rho ,B,b,A)\in ${\bf Aux}.  Suppose
that $B\in \catgen$ and that the
upper ramification set is negligible.  Also, suppose that if
$x\in B$ and $g\in G$, then either $g\cdot x=x$ or $x$ and
$g\cdot x$ can be separated by open sets.  Then $b$ is a
pseudo\'etale \catgen -morphism.
\end{theorem}

Proof will require several lemmas.

Let us begin with some point-set topology.  Let
$(G,\rho ,B,b,A)\in ${\bf Aux}.  Let $K$ be the upper ramification
set, and assume that $B-K\in \catlc$.  Put
$I=b(K)$.  Suppose
$V$ is an open connected subset of $A-I$, and let $W$ be a connected
component of $b^{-1}V$.  By elementary methods, it follows that
\begin{pflist}
\item	$b(W)=V$, and
\item	$b^{-1}V=\bigcup _{g\in G} g\cdot W$.
\end{pflist}
We shall use this observation repeated.

Let $(G,\rho ,B,b,A)\in ${\bf Aux}.  Suppose that $B\in \catlc$ and
that the upper ramification set $K$ is negligible in $B$.  Since $b$
is open, $A\in \catlc$.  Put
$I=b(K)$.  Let
$V$ be an open, conected, non-empty subset of $A$, and let $W$ be
a connected component of $b^{-1}V$.  Then $W-K$ is a (non-empty)
connected component of
$(b^{-1}V)-K=b^{-1}(V-I)$, and all other components are images of
$W-K$ under action by $G$.  Consequently, $V-I$ is connected and
non-empty.  Thus, $I$ is negligible in $A$.

Let us continue under the hypothesis of the above paragraph.  First,
we may deduce that $b$ is diffuse.  In fact, a subset element of
$A$ is negligible if and only if its inverse under $b$ is negligible in
$B$.  Thus, if $f:A\longrightarrow X$ is a continuous function such
that $f\circ b$ is diffuse, then $f$ must be diffuse as well.  In other
words, $b$ has the quotient property with respect to the category
\catgen~as well as to the category \cattop !

Suppose $(G,\rho ,X,b,Y)\in ${\bf Aux} such that, for $x\in X$ and
$g\in G$ such that $g\cdot x\neq x$, it is possible to separate $x$
and $g\cdot x$ by open sets.  We adopt some notation.  Let $x\in X$.
We denote the $G$-stabilizer of $x$ by $G(x)$.  By a symmetric
neighborhood of $x$, we mean an open, connected neighborhood $U$ such that
$$
	g\cdot U=U \mbox{~for~}g\in G(x)
	\mbox{~~ and~~}h\cdot U\cap U=\emptyset
	\mbox{~for~}h\in G-G(x).
	$$
Clearly, the hypothesis that points in $x$'s orbit can be separated
implies that there is a basis at $x$ consisting of symmetric
neighborhoods.  For $x\in X$, let $X(x)$ denote the set of symmetric
neighborhoods of $x$.

Let {\bf Test} be the class of members of {\bf Aux} whose upper
ramification set is negligible and in which the orbit of any
point can be separated by open subsets.  To proceed, we must
show that {\bf Test} satisfies
\cite[(59)]{O3}.  
We have verified that each member has the basic quotient property.
Next, we consider pullbacks.

Let $(G,\rho ,B,b,A)\in${\bf Test}, and let
$c:C\longrightarrow A$ be a \catgen -morphism.  Let
$(X;\pi _{C},\pi _{B})$ be pullback $c^{-1}(B,b)$ in the
category \cattop , and let
$(\Lambda ,\cate _{1},\lambda )$ be the construction in
Corollary~\ref{xcol1}.  As a pullback in \cattop , $X$ supports a
canonical action by $G$; as a pullback in \catgen , $\Lambda$ supports
a canonical action.  The action by $g\in G$ on $\Lambda$ can be
characterized as the unique \catgen -morphism $f$ such that
$\lambda (f(x))=g\cdot \lambda (x)$ for all $x\in \Lambda$.
Tautologically, the upper ramification set of $\Lambda$ will be
a closed subset of the pullback of the upper ramification set of $B$.
Hence, the former is negligible.  To proceed, we need one
more fact: that $\pi _{C}\circ \lambda$ is a quotient map for the
action on $\Lambda$.

Let $K$ be the upper ramification set for the action of $G$ on $X$.
Let $d=\pi _{C}$, and let
$J=d(K)$.  It is known that $X-K\in \catlc$ and that $J$ is closed.
By inspection, $d$ is a quotient for action of $G$ on
$X$.  Moreover, the lower ramification set for $d$ lies in the pullback of
the lower ramification set for $b$; therefore, both are negligible.  In
addition, by inspection, members of a $G$-orbit in $X$ can be separated
by open sets.

The behavior of $\lambda$ away from $K$ is trivial.  To finish the
present step, we must show
\begin{pflist}
\label{xlist1}
\item	the function $d\circ \lambda$ is locally open at each
	$y\in \lambda ^{-1}K$, and
\item	for each $x\in K$, $\lambda ^{-1}\{ x\}$ is a non-empty
	set of elements in the same $G$-orbit.
\end{pflist}
First, we describe the action by $G$ on $\Lambda$ explicitly.
For $g\in G$ and $(x,f)\in \Lambda$, define $g\cdot (x,f)$ to be
$(g\cdot x,h)$ where $h$ is defined by $h(U,I)=g\cdot f(U,I)$.  By
inspection, for a fixed $g\in G$, $(x,f)\mapsto g\cdot (x,f)$ is a
continuous function which satisfies the necessary commutation.

For the moment, fix $x\in X$. Suppose $U\in X(x)$.  Now $d(U)$
is connected and open, as must be $d(U)-J$.  Let $M$ be a connected
component of $d^{-1}(d(U)-J)$.  All components are $G$-conjugates
of $M$.  The union of elements $g\cdot \overline{M}$ over
$G$ is a $G$-invariant closed set, and so its image is a closed set.
That image contains all of $d(U)$.  Consequently, for any $y\in U$,
there is a component whose closure contains $y$.  Use
$y=x$, and, without loss of generality, we may assume
$x\in \overline{M}$.  Let $\Omega (U)$ be the set of connected
components of $d^{-1}(d(U)-J)$ whose closures contain $x$.

Since $U$ is symmetric, we know that $U$ itself is a connected component of
$d^{-1}d(U)$.  Consequently, every member of $\Omega (U)$ is
contained in $U$.  Now for $g\in G$, we see that either $g\in G(x)$,
in which case $g$ permutes the members of $\Omega (U)$, or
$g\notin G(x)$, in which case $g\cdot M$ is disjoint from $U$.
It follows that $U-K$ is exactly the disjoint union of the
members of $\Omega (U)$.

For each $U\in X(x)$, the size of $|\Omega (U)|$ is bounded above by
$G(x)$.  Therefore, we may choose $U_{0}\in X(x)$ for which
$|\Omega (U_{0})|$ is maximal.  That is,
\begin{equation}
	|\Omega (U_{0})|=\mbox{max}\{ |\Omega (U)|~:~U\in X(x)\}
\end{equation}
Then, for each $V\in X(x)$ such that $V\subseteq U_{0}$, it is
easily checked that
\begin{onecond}
\item	The rule $M\mapsto M\cap V$ defines a  bijection
	$\Omega (U_{0})\longrightarrow \Omega (V)$.
\end{onecond}
Once this claim is accepted, it is easily argued that
\begin{onecond}
\item	For each $M\in \Omega(U_{0})$, there is a unique
	member $f_{M}\in \Lambda (x)$ with the property that 
	$f_{M}(U_{0},U_{0}\cap K)=M$.
\end{onecond}
An immediate corollary is that the function 
$M\mapsto (x,f_{M})$ is $G(x)$-equivariant.  Condition
(\ref{xlist1}.b) follows.

Let us draw one more conclusion.
\begin{onecond}
\label{wexcomp}
\item	For $x\in X$, $U\in X(x)$, and for $M$ a connected
	component of $U-K$, there is some $f\in \Lambda (x)$
	such that $f(U,U\cap K)=M$.
\end{onecond}
Proof is left to the reader.

It remains to show that $d\circ \lambda$ is locally open.  Let
$(x,f)\in \Lambda$, and let $U\in X(x)$.  It suffices to show that
$d\circ \lambda (N(U,x,f))$ is open in $C$.  The image
$d(U)$ is known to be open.  Thus, it suffices to show that
$d\circ \lambda (N(U,x,f))=d(U)$.  That is, for each $y\in U$, we must
find an $h\in \Lambda (y)$ and $g\in G$ such that
$g\cdot (y,h)\in N(U,x,f)$.

Suppose $y\in U$.  We have already observed that every member of
$U$ lies in the closure of a $G(x)$-conjugate of
$M=f(U,U\cap K)$.  We are free to replace $y$ by any $G$-conjugate;
hence, assume that $y$ lies in the closure of $M$.

Let
$V\in X(y)$ be a symmetric neighborhood such that
$V\subseteq U$.  By (\ref{wexcomp}), there is $h\in \Lambda (y)$
such that $h(U,U\cap K)=M$.  It is clear that $h(U,I)=f(U,I)$ for every
$I\in \cate [U]$.  Thus, $(y,h)\in N(U,x,f)$.

At this point, we may legitimately state that
\begin{onecond}
\item	The class {\bf Test} satisfies
	\cite[(59.A,B,C)]{O3}.  
\end{onecond}
Let us turn to condition (D).

Let $(G,\rho ,B,b,A)\in${\bf Test}.  Let $Y\in \catgen$ be non-empty
and connected, and let $\alpha ,\beta :Y\longrightarrow B$ be
two \catgen -morphisms such that $b\circ \alpha =b\circ \beta$.
Let $K$ be the upper ramification set of $B$.  Then the subset
$Y_{1}=Y-(\alpha ^{-1}K)-(\beta ^{-1}K)$ is non-empty, connected,
and dense in $Y$.  Choose and $y\in Y_{1}$.  Then there is a unique
$g\in G$ such that $g\cdot \alpha (y)=\beta (y)$.  As the
restriction of $b$ to $B-K$ is a local homeomorphism, and since
$Y_{1}$ is connected, it follows that $g\cdot \alpha (y)=\beta (y)$
for all $y\in Y_{1}$.  Now, for $y\notin Y_{1}$, the hypothesis that
two distinct members of $\alpha (y)$'s orbit can be separated by
open sets will rule out the possibility that
$g\cdot \alpha (y)\neq \beta (y)$.  (This elementary step uses the
fact that $Y_{1}$ is dense.)  This finishs
\cite[(59.D)]{O3}.  

We now have a family of quotients which satisfies the hypothesis of
\cite[Lemma~46]{O3}.  
Let $(G,\rho ,B,b,A)\in${\bf Test}.  Let $u:U\longrightarrow B$
be an open embedding.
Corollary~\ref{xcol4} states that $b$ is open in the
traditional sense.  Let
$w:W\longrightarrow A$ be an open embedding whose image, in
the usual sense, is the image of $b\circ u$.  Then a pullback of this
$w$ will be an open embedding whose image is the set of all
$x\in B$ which are $G$-conjugate to something in the image of
$u$.  Consequently, this choice of $w$ meets condition
 \cite[(60)]{O3}.  
Therefore, by
\cite[Corollary~47]{O3}, 
every member of {\bf Test} is a perfect quotient.

At this point, Theorem~\ref{xtheo2} is a special case of
\cite[Theorem~51]{O3}.  

\section{Lifting Faithful Functors}
\label{faith}
Let \catm~be a category of orbifolds which comes with a
topological model.  In our language, this means there is a chosen
continuous functor
$\Gamma :\catm \longrightarrow \catgen [S]$ for some set $S$.
The functor $\Gamma$ is rarely full.  That is, if $X,Y$ are
topological spaces identified with orbifolds, then one does not
expect that every continuous
$f:X\longrightarrow Y$ will lift to a morphism of orbifolds.
However, in many situations, $\Gamma$ is expected to be faithful.
That is, if $f,g:M\longrightarrow N$ are two morphisms of orbifolds
such that $\Gamma (f)=\Gamma (g)$, then
$f=g$.  Informally, we say $f$ and $g$ are equal if their ``underlying''
functions agree.

In our present program, a useful category \catm~appears as
the expansion of an intial category \catc ---that is,
\catm~is \catcpl ,
${\cal C}^{++}$, or some higher iterate of the plus construction.
On the initial
\catc , it should be clear, by inspection, that $\Gamma$ is faithful.
In this section, we prove that all the lifts of $\Gamma$ are faithful,
under a reasonable hypothesis.

Our standing hypothesis is
\begin{onecond}
\label{xhyp3}
\item	Let \catc~be a topologized category which satisfies
	\cite[(11)]{O3}. 
	 Let $+:\catc \longrightarrow \catcpl$ be a plus
	 construction, and let $p$ and $s$ denote the component
	 pasting and smoothing functors, respectively.  Regard
	 $\catc ^{p}$ as topologized in the usual way.
\end{onecond}
Let us begin with an elementary reduction.

\begin{lemma} \label{xlem6}
Assume (\ref{xhyp3}).  Let \cate~be a topologized category and let
$\Gamma :\catcpl \longrightarrow \cate$ be a covariant functor.
Assume that
\cate~is quasi-intrinsic, in the sense of \cite[Definition 12.5]{LG},
and that $\Gamma \circ s$ is continuous.  Then $\Gamma$ is
faithful if and only if for each $B\in \catc$ and each
$A_{0}\in \catcp$, the function
$$
	\morcp{B^p}{A_0}\longrightarrow
	\mbox{Mor}_{\cal E}(\Gamma (B^{+}),\Gamma (A_{0}))
	$$
defined by $\Gamma$, is injective.
\end{lemma}

\begin{pfproof}
Let $A_{0},B_{0}\in $Obj(\catcpl )=Obj(\catcp ).  Let
$f_{1},g_{1}:B_{0}\longrightarrow A_{0}$ be \catcpl -morphisms.
Then there exist \catcp -morphisms
$$\begin{array}{rrr}
	x_{0}:X_{0}\longrightarrow B_{0},  &&
	 f_{0}:X_{0}\longrightarrow A_{0},   \\
	y_{0}:Y_{0}\longrightarrow B_{0}~   &
	\mbox{and} &  g_{0}:Y_{0}\longrightarrow A_{0},
	\end{array}
	$$
such that $x_{0}$ and $y_{0}$ are pseudoisomorphisms and
$$
	f_{0}^{s}=f_{1}\circ x_{0}^{s}
	\mbox{~~and~~} g_{0}=g_{1}\circ y_{0}^{s}.
	$$
Let $((P,p_{0});r_{0},t_{0})$ be a fibered product
$(X_{0},x_{0})\times _{B_{0}}(Y_{0},y_{0})$ in $\catc ^{p}$.  Then
$p_{0}$, $r_{0}$ and $t_{0}$ are all pseudoisomorphisms and
$$
	\{ f_{0}\circ r_{0}\} ^{s}=f_{1}\circ p_{0}^{s}
	\mbox{~~ and~~}
	\{ g_{0}\circ t_{0}\} ^{s}=g_{1}\circ p_{0}^{s}.
	$$
Smoothing takes $p_{0}$, $r_{0}$ and $t_{0}$ to isomorphisms.  Thus,
$\Gamma (f_{1})=\Gamma (g_{1})$ if and only if
$\{ \Gamma \circ s\} (f_{0}\circ r_{0})=
\{ \Gamma \circ s\} (g_{0}\circ t_{0})$.

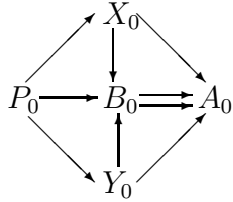
\begin{figure}
\caption{Plus Morphisms as Maps between Canopies}
\label{afig1}
\centering
\begin{picture}(90,100)(-4,-36)
\put(-2,0){$P_{0}$}
\put(34,32){$X_{0}$}
\put(34,0){$B_{0}$}
\put(70,0){$A_{0}$}
\put(34,-32){$Y_{0}$}

\put(5,10){\vector(1,1){26}}
\put(6,-3){\vector(1,-1){25}}
\put(10,4){\vector(1,0){22}}
\put(48,1){\vector(1,0){20}}
\put(48,5){\vector(1,0){20}}
\put(38,30){\vector(0,-1){20}}
\put(40,-22){\vector(0,1){20}}
\put(47,36){\vector(1,-1){26}}
\put(47,-28){\vector(1,1){26}}

\end{picture}
\end{figure}
It follows that $\Gamma$ is faithful if and only if
$\Gamma \circ s$ is faithful.

Let $A_{0},B_{0}\in \catcp$ and let
$f_{0},g_{0}\in $Mor$_{C^{p}}(B_{0},A_{0})$.  If
$f_{0}\circ \iota _{j}=g_{0}\circ \iota _{j}$ for each
$j\in \Lambda (B_{0})$, then $f_{0}=g_{0}$.  Now
$\Gamma \circ s$ is continuous, so it sends the assigned cover to
$B_{0}$ to an \cate -cover.  Since the topology of
\cate~is quasi-intrinsic, it is also true that if
$\Gamma (f_{0})\circ \Gamma (\iota _{j})=
\Gamma (f_{0}\circ \iota _{j})$ and
$\Gamma (g_{0})\circ \Gamma (\iota _{j})=
\Gamma (g_{0}\circ \iota _{j})$ agree for every $j$, then
$\Gamma (f_{0})=\Gamma (g_{0})$.  The desired conclusion follows.
\end{pfproof}

For convenience, we introduce a term for use in this section only.
Let \catd~be a topologized category.  We say \catd~meets the
{\em m-hypothesis\/} if the following two conditions are true.
\begin{pflist}
\label{xlist2}
\item	Every formal \catd -subset $b$ can be decomposed as
	$b=w\circ f$ where $w$ is a monomorphic formal
	subset and $f$ overlays absolutely.  (Recall that, in
	this situation, $f$ is a pullback for $b$.)
\item	Let $\theta$ be an indexed cone of formal
	\catd -subsets into an object $X$.  Then the
	canopy of $\theta$ admits an affinization.
\end{pflist}
Condition (\ref{xlist2}.b) is \cite[(9.c)]{O3}, 
albeit in a more general context.

\begin{proposition}
\label{xprop3}  Assume (\ref{xhyp3}).

(A)  The category \catcpl , with the e.l.-topology, meets the m-hypothesis.

(B)  Suppose that \catc~meets the m-hypothesis.  Let $A\in \catc$,
and let
$b_{1}:B_{0}\longrightarrow A^{+}$ be a monomorphic e.l.-subset.
Then there is a monomorphic formal \catc -subset
$b:B\longrightarrow A$ such that $b_{1}$ is
$\catcpl /A^{+}$-isomorphic to $b^{+}$.
\end{proposition}

\begin{pfproof}
We sketch the proof.  Details are in the style of \cite{LG}.

We begin with (\ref{xlist2}.a).  Let
$b_{1}:B_{0}\longrightarrow A_{0}$ be an e.l.-subset.  We are
free to replace $b_{1}$ by any morphism
$\catcpl /A^{+}$-isomorphic to it.  Thus, we may assume that
$b_{1}=b_{0}^{s}$ where $b_{0}:B_{0}\longrightarrow A_{0}$ is a
\catcp -morphism such that, for each $j\in \Lambda (B_{0})$,
$b_{0}\circ \iota _{j}$ is a formal ${\cal C}^{p}$-subset.  Now
a product of two affine formal subsets is known to be affine.
It follows that there is a
${\cal C}^{p}$-morphism $w_{0}:W_{0}\longrightarrow A_{0}$ such that
\begin{pflist}
\item	$\Lambda (W_{0})=\Lambda (B_{0})$,
\item	for each $j\in \Lambda (W_{0})$,
	$w_{0}\circ \iota _{j}=b_{0}\circ \iota _{j}$,~~~and
\item	for $(j,k)\in \Lambda (B_{0})^{2}$, the image of
	$(W_{0}[j,k];\rho _{1},\rho _{2})$ under
	pasting is a product
	$\{ b_{0}\circ \iota _{j}\} \times _{A_{0}}
	\{ b_{0}\circ \iota _{k} \}$.
\end{pflist}
It follows that
\begin{pflist}
\item	$w_{0}^{s}$ is a monomorphic e.l.-subset,
\item	there is a morphism $f_{0}:B_{0}\longrightarrow W_{0}$
	such that, for each $j\in \Lambda (W_{0})$,
	$f_{0}\circ \iota _{j}=\iota _{j}$,
\item	$f_{0}^{s}$ is an e.l.-subset which covers.
\end{pflist}
Thus, $b_{0}^{s}=w_{0}^{s}\circ f_{0}^{s}$ is a decomposition
of the required type.

Let us interpolate a proof of Part (B).  In the
present context, suppose that
\catc~meets the m-hypothesis and that $A_{0}=A^{p}$ where
$A\in \catc$.  Then $W_{0}$ is a canopy of a cone of formal
\catc -subset.  By assumption, it admits an affinization in
\catc .  Let $w$ be the affinization of $w_{0}$.  Then $w$ is a
monomorphic local subset.  By condition
\cite[11.D]{O3}, 
$w$ is a local subset.

Now return to Part (A).  Suppose $A_{0}\in \catcpl$ and $\theta$ is
a cone of
e.l.-subsets into $A_{0}$ indexed by a set $T$.  Again, we may assume
without loss of generality that for each $t\in T$,
$\theta (t)=^{t}b_{0}^{s}$ where
$^{t}b_{0}:^{t}B_{0}\longrightarrow A_{0}$ is a \catcp -morphism
such that $^{t}b_{0}\circ \iota _{j}$ is a formal \catcp -subset for
each $j\in \Lambda (^{t}B_{0})$.  Put
$$
	J=\{ (t,j)~:~t\in T,~j\in \Lambda (^{t}B_{0})\} .
	$$
There is $U_{0}\in \catcp$ of type Int$(J)$ and a \catcp -morphism
$u_{0}:U_{0}\longrightarrow A_{0}$ for which
\begin{pflist}
\item	for each $(t,j)\in J$,
	$u_{0}\circ \iota _{(t,j)}=^{t}b_{0}\circ \iota _{j}$,
\item	for each pair $((t,j),(r,k))\in J^{2}$, the image of
	$(U_{0}[(t,j),(r,k)];\rho _{1},\rho _{2})$ under
	the pasting functor is
	a fibered product $\{ ^{t}b_{0}\circ \iota _{j}\}
	\times _{A_{0}}\{ ^{r}b_{0}\circ \iota _{k}\}$.
\end{pflist}
For each $t\in T$, there is a unique \catcp -morphism
$^{t}p_{0}:^{t}B_{0}\longrightarrow U_{0}$ such that, for each
$j\in \Lambda (^{t}B_{0})$, $^{t}p_{0}\circ \iota _{j}=\iota _{(t,j)}$.
Then $t\mapsto ^{t}p_{0}^{s}$ is an affinization for $\theta$ in \catcpl .
\end{pfproof}

Let us recall some facts and comments from \cite{O2}.  Let \catd~be
a topologically componentwise category whose topology is intrinsic
and flush.  In \catd , a morphism whose domain is connected is said
to be ``connected''; similarly, it is called non-empty if its domain
is ``non-empty''.  A
{\em decomposition into connected component\/} of an object $A$
is an indexed family $\theta$ of complemented morphisms into
$A$ for which
\begin{pflist}
\item	$\theta$ is a cover, and
\item	for $j,k\in \dom{\theta}$ such that $j\neq k$, the product
	$\theta (j)\times _{A}\theta (k)$ is empty,
\item	for each $j\in \dom{\theta}$, the domain of
$\theta (j)$ is connected and non-empty.
\end{pflist}
Recall from \cite{O2} that, in such a category, every object with a
cover by connected morphisms admits a decomposition into connected
components.

\begin{lemma}
\label{xlem7} Let \catc~and \cate~be topologized categories.
Assume both are intrinsic and topologically componentwise.  By a
{\em topologically componentwise\/} functor from
$\catc \longrightarrow \cate$, we mean a continuous functor
$\Gamma$ from \catc~to \cate~such that
\begin{pflist}
\item	$\Gamma$ sends empty \catc -objects to empty \cate -objects,
\item	$\Gamma$ sends non-empty connected \catc -objects to
non-empty connected \cate -objects.
\end{pflist}
Let $\Gamma :\catc \longrightarrow \cate$ be a topologically
componentwise functor.  Assume that every
\catc -object has a cover by connected morphisms.

(A)  $\Gamma$ sends complemented \catc -morphisms to
complemented \cate -morphisms.

(B)  $\Gamma$ sends non-empty objects to non-empty objects.

(C)  Let $\theta$ be a decomposition into connected
components of an object $A\in \catc$.  Then
$\Gamma \circ \theta$ is a decomposition into connected
components of $\Gamma (A)$.
\end{lemma}

\begin{pfproof}
Every \catc -object has a decomposition into connected components.
Thus, every non-empty
\catc -object is a the codomain of a connected morphism.  It
follows that $\Gamma$ sends non-empty objects to non-empty objects.

In either \catc~or \cate , a complemented morphism can be
characterized as a morphism
$b:B\longrightarrow A$ such that
\begin{pflist}
\item	$b$ is a formal subset,
\item	$(B;1_{B},1_{B})$ is a self-product $b\times _{A}b$,
\item	there is another morphism
	$c:C\longrightarrow A$ which satisfies (a) and
	(b) and for which
	$b\times _{A}c$ is empty and $\{ b,c\}$ is a cover.
\end{pflist}
These properties are preserved by $\Gamma$.  Part (A) follows.
Part (C) is an easy consequence.
\end{pfproof}

\begin{theorem} \label{big1}
Assume (\ref{xhyp3}).  Let \cate~be a topologized category, and let
$\Gamma :\catcpl \longrightarrow \cate$ be a
covariant functor.  Assume
\begin{pflist}
\label{xlist3}
\item	\catc~is topologically componentwise and every
	\catc -object has a cover by connected
	\catc -morphisms,
\item	every formal \catc -subset is discrete,
\item	\catc~satisfies the m-hypothesis,
\item	the topology of \cate~is intrinsic, flush and
	topologically componentwise,
\item	every formal \cate -subset is discrete,
\item	$\Gamma \circ +$ is faithful and topologically componentwise,
\item	$\Gamma \circ s$ is continuous,
\item	for $b$ a formal \catc -subset, if
	$\Gamma (b^{+})$ is an \cate -isomorphism, then
	$b$ is a \catc -isomorphism.
\end{pflist}
Then
\begin{pflist}
\label{xlist4}
\item	$\Gamma$ is a faithful functor, and it sends
	connected \catcpl -objects
	to connected \cate -objects,
\item	if $b:B\longrightarrow A$ is a formal \catc -subset and
	$c:C\longrightarrow A$ is an arbitrary \catc -morphism with
	connected domain, then the function induced by $\Gamma$
	$$
	\mbox{Mor}_{C/A}(C,c;B,b)\longrightarrow
	\mbox{Mor}_{E/\Gamma (A)}(\Gamma (C,c),\Gamma (B,b))
	$$
	is a bijection,
\item	if $b$ is an e.l.-subset such that $\Gamma (b)$ is an
	\cate -isomorphism, then $b$ is a \catcpl -isomorphism.
\end{pflist}
\end{theorem}

\begin{remark}
In practice, attention is limited to topologies meeting
an explicit list of axioms.  The theorem applies to topologies
which are flush, intrinsic, topologically componentwise and such
that every formal subset is discrete.  In a practical situation, there
would be a canonical choice of topology for \catcpl , related to the
e.l.-topology, with respect to which $\Gamma$ is continuous.
In this context, conditions (\ref{xlist3}.a,b,c,d,e,g) are true by
fiat or by elementary arguments (eg., using the fact that the
e.l.-topology meets the m-hypothesis to show that the chosen
topology of \catcpl~meets the same condition).  The ``real''
hypothesis is (\ref{xlist3}.f,h).  The theorem says that if the
hypothesis holds for
\catc , it will hold for \catcpl .
\end{remark}

\begin{pfproof}
The conclusions of Lemma~\ref{xlem7} apply to
$\Gamma$ by (\ref{xlist3}.g).

Let $A_{0}$ be a connected \catcpl -object.  For each
$j\in \Lambda (A_{0})$, let $\phi _{j}$ be a decomposition of
$A_{0}[j]$ into connected components.  Let $\lambda$ be the
refinement of the assigned cover through $j\mapsto \phi _{j}$,
and let $\Omega$ be the domain of $\lambda$.  Note that no
member of $\lambda$ is empty.  On $\Omega$, define the
linking relation, as in \cite{O2}.  That is, define
$\sim$ on $\Omega$ to be the smallest equivalence relation such that
$r\sim s$ when $\lambda (r)\times _{A_{0}}\lambda (s)$
is non-empty.

We now cite results from \cite{O2}.  Since $A_{0}$ is connected in
\catcpl , all members of $\Omega$ are equivalent under
$\sim$.  By Lemma~\ref{xlem7}, the image of $\lambda$ under
$\Gamma \circ s$ is a connected cover of $\Gamma (A_{0})$.
Define $\approx$ on $\Omega$ to be the smallest
equivalence relation such that
$r\approx s$ when
$\Gamma (\lambda (r))\times _{\Gamma (A_{0})}
\Gamma (\lambda (s))$ is non-empty.
Recall that $\Gamma \circ s$ is continuous by assumption.
It follows that $\approx$ is $\sim$, and
that all members of $\Omega$ are equivalent.  Therefore,
$\Gamma (A_{0})$ is connected.

Let $B\in \catc$, $A_{0}\in \catcp$ and
$f_{0},g_{0}\in $Mor$_{C^{p}}(B^{p},A_{0})$ such that
$\Gamma (f_{0})=\Gamma (g_{0})$.  Take
$j,k\in \Lambda (A_{0})$ and \catc -morphisms
$f:B\longrightarrow A_{0}[j]$ and $g:B\longrightarrow A_{0}[k]$
such that $f_{0}=\iota _{j}\circ f^{p}$ and
$g_{0}=\iota _{k}\circ g^{p}$.  Since $\Gamma \circ s$ is a
continuous functor, $(\Gamma (A_{0}[j,k]);
\Gamma (\rho _{1}),\Gamma (\rho _{2}))$ is a fibered product
$\Gamma (\iota _{j})\times _{\Gamma (A_{0})}
\Gamma (\iota _{k})$.  Since
$\Gamma (g_{0})=\Gamma (\iota _{k})\circ \Gamma (g^{+})$ and
$\Gamma (f_{0})=\Gamma (\iota _{j})\circ \Gamma (f^{+})$
agree, there is an \cate -morphism
$H:\Gamma (B)\longrightarrow \Gamma (A_{0}[j,k])$
such that $\Gamma (\rho _{1})\circ H=\Gamma (f^{+})$ and
$\Gamma (\rho _{2})\circ H=\Gamma (g^{+})$.  By
Lemma~\ref{xlem6}, to prove that $\Gamma$ is faithful it
suffices to show that, in this situation,
$H=\Gamma (h)$ for some
\catc -morphism $h:B\longrightarrow A_{0}[j,k]$.  Thus,
(\ref{xlist4}.a) will be true if (\ref{xlist4}.b) holds.

Let us turn to (\ref{xlist4}.b).  We begin with a comment
equally applicable to \catc~and \cate.

Let \catd~be any componentwise category.  Let
$b:B\longrightarrow A$ be a discrete pullback base, let
$C$ be a connected object, and let
$c:C\longrightarrow A$ be a \catd -morphism.  Let $(P;p,q)$
be a pullback
$c^{-1}(B,b)$.  It is known, from \cite{O1}, that there is a
bijection from the set $S$, of
$\catd /A$-morphisms $f:(C,c)\longrightarrow (B,b)$, to the
set $T$, of connected components of $P$ on which $p$ is
an isomorphism, which assigns to each
$f\in S$ the connected component of $P$ through which
$1_{C}\times f$ factors.

Suppose
$b:B\longrightarrow A$ is a formal \catc -subset,
$c:C\longrightarrow A$ is a connected morphism and
$(P;p,q)$ is a pullback
$c^{-1}(B,b)$.  The functor $\Gamma$ preserves these
relations and sends components of $P$ to components of
its image.  The map in (\ref{xlist4}.b) can fail to be
bijective only if there is a connected component of
$x:X\longrightarrow P$ such that $p\circ x$ is not a
\catc -isomorphism but $\Gamma (p\circ x)$ is an
\cate -isomorphism.  The latter is explicitly disallowed
by assumption (\ref{xlist3}.h).

Conclusions (\ref{xlist4}.a,b) have been verified.  The last is
comparatively simple.  Suppose
$b_{1}:B_{0}\longrightarrow A_{0}$ is an e.l.-morphism for
which $\Gamma (b_{1})$ is an isomorphism.  Let $\theta$ be
an affine
\catcp -cover of $A_{0}$.  To show that $b_{1}$ is an
isomorphism, it suffices to prove that each pullback along a
member of $\theta$ is an isomorphism.  Hence, we may assume
that $A_{0}=A^{+}$ for some $A\in \catc$.

It is known that $\Gamma$ is faithful.  Thus, if $b_{1}$ is
not monomorphic in \catcpl , its image cannot be monomorphic.
Therefore, $b_{1}$ is a monomorphic e.l.-subset into an affine
object.  From Proposition~\ref{xprop3}, it follows that $b_{1}$ is
isomorphic to a formal \catc -subset.  By hypothesis (\ref{xlist3}.g),
$b_{1}$ is an isomorphism.
\end{pfproof}

\section{What is a Morphism?}
\label{whatis}
Let $M$ and $N$ be two $C^{\infty}$-manifolds.  Let
$(G,\sigma )$ and $(H,\tau )$ be finite group actions on
$M$ and $N$ respectively, which are discrete with respect to diffuse
$C^{\infty}$-functions between manifolds.  This Section
discusses the specific issue of what is a morphism between
the quotient spaces $G\backslash M\longrightarrow H\backslash N$.
For us, this amounts to looking at the abstract
machinery and rewriting the formal definition in this explicit
case.  Irregularities and surprising behaviors
have been noted in explicit situations.  We offer our spin on
a pathological situation discussed by Schwarz in
\cite{GS}.

To begin, let us lay out the above situation in the abstract.  Let
\catm~be a category which is being used to generate orbifolds.  That is,
\begin{pflist}
\item	\catm~is a topologized category,
\item	the topology of \catm~is flush, intrinsic and
	topologically componentwise,
\item	every \catm -object has a cover by connected objects,
\item	the topology of \catm~is {\em pseudo\'etale\/} (or
	{\em torsorial\/}); that is, a local \catm -subset
	which is discrete and finite (respectively, just discrete)
	must be a formal \catm -subset),
\item	\catm~meets the set-theoretic axioms in
	\cite[(11)]{O3}. 
\end{pflist}
By assigning to $\catm ^{+}$ the class of discrete (and, depending
on context, finite)
e.l.-subsets, we arrange for $\catm ^{+}$ to meet the same axioms.
Thus, the fundamental construction generates a sequence of
enlargements \catm ,
$\catm ^{+}$, $(\catm ^{+})^{+}$,... whose members are orbifolds.
Assume $S$ is a set and
$\Gamma :\catm \longrightarrow \catgen [S]$ is a faithful
continuous functor.  We also assume that \catm~meets the
m-hypothesis, and that, for $b$ a
formal \catm -subset, $b$ is an \catm -isomorphism if and only if
$\Gamma (b)$ is a $\catgen [S]$-isomorphism.  Then
Theorem~\ref{big1} applies to all lifts of $\Gamma$.

 In terms of the particular example, \catm~is the category whose
 objects were all $C^{\infty}$-manifolds but whose morphisms
 consisted only of diffuse
$C^{\infty}$-functions, formal \catm -subsets are finite
covering maps whose Jacobian is invertible at each
point, $S=\bbr$, and $\Gamma$ is the forgetful
functor.

Let $M,N\in \catm$, and let $(G,\sigma )$ and $(H,\tau )$ be
group actions on $M$ and $N$, respectively, which are discrete
in the categorical sense.  Let $M_{0}$ be the canopy of Type
Int($\{ 1\}$) such that $M_{0}[1]=M$, $M_{0}[1,1]$ is the disjoint
union of copies of \catm , indexed by $G$, and, for $g\in G$, the
restrictions of $\rho _{1}$ and $\rho _{2}$ to the $g$-th copy of
$M$ are $1_{M}$ and
$\sigma (g)$, respectively.  Let $N_{0}$ denote the analogous
canopy for $N$ and $H$.  Then, as members of $\catm ^{+}$,
$M_{0}$ is the quotient $G\backslash M$ and $N_{0}$ is the
quotient $H\backslash N$.  Let $\Gamma ^{+}$ denote the
continuous extension of $\Gamma$ to
$\catm ^{+}\longrightarrow \catgen [S]^{+}$.

Before we continue with the abstract situation, let us look at a particular
choice of manifolds and actions.  Schwarz, in \cite{GS}, illustrates
surprising behavior with a morphism from
$\{ -1,1\} \backslash \bbr ^{2}
\longrightarrow \{ -1,1\} \bbr$.  As noted in
Subsection~\ref{nonquot}, the canonical action by $\{ -1,1\}$ on
\bbr~is not discrete in our category.  However, we get the same effect
by taking his situation and forming products with \bbr .  Also, we
have to worry that our functions are diffuse,
a concept which is not in Schwarz's work.
Below is the update of his example.

Let $M=\bbr ^{3}$, $N=\bbr ^{2}$, $H=G=\{ -1,1\}$, and
characterize the two actions by
$$
	-1\cdot (x,y,z)=(-x,-y,z)
	\mbox{~~ and~~} -1\cdot (x,z)=(-x,z).
	$$
Let $f:(0,\infty )\longrightarrow \bbr$ be a
$C^{\infty}$-function such that
\begin{pflist}
\item	$f$ is periodic of period 4,
\item	$f^{\prime}$ is negative on intervals (1,2) and (2,3),
	and is positive on (3,4) and (4,5),
\item	all derivatives of $f$ vanish at every even integer,
\item	$f(2k)=0$ for all $k\in \bbn$.  Equivalently, $f$ is
	positive on $(4k,4k+2)$ and negative on
	$(4k+2,4k+4)$ for all
	$k\in \bbn \cup \{ 0\}$.
\end{pflist}
For each $x\in (0,\infty )$, define
$g(x)=e^{-n}f\left( \frac{1}{x}\right)$ where $n$ is the smallest
even integer greater than or equal to $1/x$.  Adopt the convention
that $g(0)=0$.  Define $F:M\longrightarrow N$ by
\begin{equation}
	\begin{array}{r}
	F(r\cdot \mbox{cos}(\theta ),r\cdot
	\mbox{sin}(\theta ),z)~=~
		\left\{ \begin{array}{cl}
		(g(r)\cdot\mbox{sin}(\theta ),z) &
			\mbox{if~}g(r)>0 \\
		(g(r)\cdot\mbox{sin}(2\theta ),z) &
			\mbox{if~}g(r)\leq 0
		\end{array} \right. \\
	\mbox{for all~}(r,\theta ,z)\in [0,\infty )
	\times [0,2\pi ]\times \bbr .
	\end{array}
\end{equation}
Standard theory verifies that $F$ is $C^{\infty}$.

We must verify that $F$ is diffuse.  On any subset where
$F$ can be identified with a projection
$U\times V\longrightarrow V$,
$F$ is diffuse.  Also, if $F$ is diffuse on each member of a
family of open subsets, then it is diffuse on the union.  Thus,
restriction of $F$ to the set of points where its Jacobian
is surjective is diffuse.

In $M$, let $\ell$ be the line $\{ (0,0)\} \times \bbr$.  For
each $n\in \bbn$, let $C_{n}$ be the cylinder of all points
distance $1/(2n)$ from $\ell$.  Let $C$ be the union of
these cylinders with $\ell$.  For each odd integer $k$, let
$\Xi _{k}$ be union of lines of points of the form $(1/k,a,z)$
where $a\in \{ \pi /2,3\pi /2\}$ or
$a\in \{ \pi /4,3\pi /4,5\pi /4,7\pi /4\}$, depending on the
sign of $g(1/k)$.  Let
$\Xi$ be the union of these sets of lines.  A brief calculation
verifies that
$F$ is diffuse on $M$ minus $C\cup \Xi$.  In fact, our
characterization of $f$ allows us to rule out the cylinders.  Although
the Jacobian degenerates, one can make a continuous identification
with projection.  Instead, let us use this example as an excuse to
introduce two practical lemmas.

\begin{proposition}
\label{xprop4}  Let $X,Y\in \catgen$ and let
$f:X\longrightarrow Y$ be a continuous function.  Let
$I\subseteq X$ be a negligible subset.  If the restriction of
$f$ to $X-I$ is diffuse, then $f$ is diffuse.
\end{proposition}

\begin{pfproof}
This is a trivial consequence of Corollary~\ref{lcor1}.
\end{pfproof}

Being diffuse is a local property, so we may invoke
Proposition~\ref{xprop4} on each line in $\Xi$.  Therefore,
$F$ is diffuse on
$M$ minus $C$.  In fact, we can go further.  The axis
$\ell$ is negligible in $M$, so to prove that $F$ is diffuse on
$M$, it suffices to prove the diffuse property on
$M-\ell$.  Consequently, it suffices to prove that $F$ is
locally diffuse about each individual cylinder $C_{n}$.

\begin{proposition}
\label{xprop5}  Let $X$ be a topological manifold, let
$Y\in \catgen$ be normal (and Hausdorff) in the topological
sense, and let
$f:X\longrightarrow Y$ be a continuous function.  Let $Z$
be a submanifold of $X$ of codimension 1.  Suppose $f$ is
diffuse on $X-Z$.  Then $f$ is diffuse unless and only unless
there is a non-empty open subset $W$ of $Z$ such that the
closure of $f(W)$ is negligible in $Y$.
\end{proposition}

\begin{pfproof}
Suppose $f$ is not diffuse.  Let $(U,I)$ be a negligible
subset-element of $Y$ whose inverse is not negligible.  Now
$f^{-1}I-Z$ is known to be negligible in $f^{-1}U-Z$.  Thus,
$f^{-1}U\cap Z$ must not be negligible in $f^{-1}U$.  There is
$x\in Z\cap f^{-1}I$ at which the local condition
of Lemma~\ref{ylem2} fails.  Because $Z$ is of codimension 1,
it follows easily that there is a neighborhood $W_{1}$ of
$x$ in $Z$ which is entirely contained if $f^{-1}I$.  Now $Y$ is
normal and Hausdorff, so there is a closed neighborhood $T$ of
$f(x)$ in $Y$ such that $T\subseteq U$.  Let
$W_{2}=W_{1}\cap f^{-1}T$, and let $W$ be the interior of this
set with respect to the topology of $Z$.  Then $x\in W$ and the
closure of $f(W)$ in $Y$ is a closed subset of $I$.  It follows that
$W$ has the stated properties.

The converse is trivial.
\end{pfproof}

On each cylinder, $F$ is projection $(x,y,z)\mapsto (0,z)$.
The closure of the image, under $F$, for any non-empty open
subset of the cylinder will contain a line segment.  In $N$, any
line segment fails to be negligible.  Hence, $F$ is diffuse.

We are ready to consider morphisms at three levels.

Let us recall how $M_{0}$ serves as $G\backslash M$.  The quotient map is
$\iota =\iota _{1}^{s}:M^{+}\longrightarrow M_{0}$.  The map
$\iota$ is purely formal.  Its construction does not include a
subtle study of quotients.  What links the purely formal to the
practical is the lift of the functor $\Gamma$.

Consider the situation when $M$ is a manifold.  Let
$q:\Gamma (M)\longrightarrow Q$ be the standard
topological quotient.  (Here, $\Gamma (M)$ is just
$M$ with differential structure omitted.)  If the fixed
point set of each non-trivial member of $G$ is negligible,
then we know that $q$ describes a pseudo\'etale quotient
with respect to \catgen .  In other words, we may choose
$\Gamma ^{+}(\iota )$ to be $q$ and $\Gamma ^{+}(M_{0})=Q$.
It is justified as, up to isomorphism, the only way to
extend $\Gamma$ continuously.

Unfortunately, it is possible that $q$ is
{\em not\/} pseudo\'etale.  This is the case if $M=\bbr ^{2}$ and
$\{ -1,1\}$ acts as described earlier.  In this case,
$\Gamma ^{+}(M_{0})$ can not be assigned a
\catgen -object without violating continuity of functors.  It
is tempting to define it to be $Q$, but, as we shall show
shortly, this decision would sacrifice the valuable property
of faithfulness.  We may say that $\Gamma ^{+}(M_{0})$ exists
in a suitable expansion $\catgen [S]^{+}$ of $\catgen [S]$, and,
in that expansion, it serves as $G\backslash \Gamma (M)$.  We may also
say that the expansion supports a fundamental group and
cohomological theories, although these differ from the classical.

Actually, for our present question, $M_{0}$ is not the
troublesome quotient.  In $\catm ^{+}$, $M_{0}$ is a true
quotient for the action by $G$.  A morphism on $M_{0}$
may be effectively defined as a $G$-invariant morphism on
$M$!  So, definition of a morphism
$M_{0}\longrightarrow N_{0}$ is really a question of
definition from
$M\longrightarrow N_{0}$.

Our study reduces to two challenges:
\begin{pflist}
\label{xlist5}
\item	describe a morphism $M\longrightarrow H\backslash N$,
\item	determine when two such morphisms agree.
\end{pflist}
Issue (\ref{xlist5}.b) runs directly into Seifert
Boundaries and faithfulness of $\Gamma$.

Formally, a morphism $M\longrightarrow N_{0}$ is represented by a pair
$(\theta ,f)$ where
\begin{pflist}
\item	$\theta$ is a cover of $M$ (in ${\cal M}$) indexed
	by some set $J$,
\item	for each $j\in J$, $f(j)$ is an ${\cal M}^{p}$-morphism
	$M\longrightarrow N_{0}$,
\end{pflist}
and certain conditions are met.  Actually, in this case,
$f(j)$ may be interpreted simply as an \catm -morphism
$M\longrightarrow N$, with the understanding that certain
morphisms are regarded as equivalent.

Even in an explicit context, such as the category of manifolds,
each map $\theta (j)$ is not limited to be just open embedding.
The family $\theta$ covers in the topology of \catm .
That topology includes finite-to-one maps.  There is a good
heuristic reason for this.  Intuitively, we wish to define a morphism
$M\longrightarrow N_{0}$, but our language only contains
morphisms between ``affine'' objects like $M$ and $N$.  A given
$f:M\longrightarrow N_{0}$ need not factor through the canonical
projection $\alpha :N\longrightarrow N_{0}$.  In this case, we can
not represent $f$ as $\alpha$ composed with something.
Instead, we try to classify $f$ by its
{\em pullback\/} along $\alpha$ to
$\pi :M^{\ast }=M\times _{N_{0}}N\longrightarrow N$.  Although
the domain of $\pi$ need not be affine, it has a cover by affines,
and $\pi$ restricted to each of these is a morphism in the old
sense.  Thus, $f$ gets expressed as a family of maps into $N$
from objects which are, crudely, open subsets of a finite cover
of $M$.  (This characterization is spiritually right,
but too simple formally.)

We do have some leeway.  We may replace any member of
$\theta$ by a refinement.  In particular, we may choose the
maps so that each has its
{\em image\/} constrained to lie in some indicated basis.  We
may pass from $\theta$ to its restriction to a subset of
$J$, provided that the restriction covers.  We may assume that
the domain of each $\theta (j)$ is connected.

If the topology is based on morphisms with non-trivial
ramification, then such morphisms may appear in
$\theta$.  The derived pseudo\'etale topology is intended to include
all forms of singularity which do not prevent certain manipulations.
In it, a cover may consist of complicated maps.  On the other hand,
the elementary finite-to-one topology is intended to be the simplest
of topologies from which quotients may be generated.
In it, every cover can be refined to a pseudogeometric cover; that is,
in any concrete case, a usual cover by open subsets.  In the
standard theory of real orbifolds, if $M$ is a manifold (rather than
an abstract orbifold) then $\theta$ may be chosen to be 
injections for all members of a cover by open subsets.

Two related questions arise.
\begin{pflist}
\label{wagree}
\item	When does a pair $(\theta ,f)$, represent
	an $\catm ^{+}$-morphism?
\item	Given pairs $(\theta ,f)$ and $(\phi ,g)$ which represent
	morphisms, when do they represent the same one?
\end{pflist}
In (\ref{wagree}.a), it is necessary and sufficient that for every
pair $(j,k)\in \dom{\theta}^{2}$ and for any $(P;\pi _{1},\pi _{2})$
a product $\theta (j)\times _{M}\theta (k)$, it is true that
$f(j)\circ \pi _{1}=f(k)\circ \pi _{2}$.  In (\ref{wagree}.b), the key
condition is that for every $j\in \dom{\theta}$, $k\in \dom{\phi}$
and any $(P;\pi _{1},\pi _{2})$
a product $\theta (j)\times _{M}\phi (k)$, it is true that
$f(j)\circ \pi _{1}=g(k)\circ \pi _{2}$.  The basic issue is
\begin{onecond}
\label{wequal}
\item	For each $U\in \catm$ and $f,g:U\longrightarrow N$
	\catm -morphisms, when do $f$ and $g$ represent the
	same morphism into $N_{0}$?
\end{onecond}
The theory gives an unambiguous answer.  However, there is an
tempting deception.

The formalism is precise.  Let $U$, $f$ and $g$ be as in
(\ref{wequal}).  Then $f$ and $g$ represent the same thing if
and only if there is an \catm -morphism
$\delta :U\longrightarrow N_{0}[1,1]$ such that
\begin{equation}
\label{weq4}
	\rho _{1}\circ \delta = f
	\mbox{~~and~~}\rho _{2}\circ \delta = g.
\end{equation}
The first composition in (\ref{weq4}) means that, on each connected
component of $U$, $\delta$ restricts to a copy of $f$ which maps
into one of the copies of $N$ in $N_{0}[1,1]$.  In other words, for
each connected component $C$ of $U$, there is $h_{C}\in H$
such that
\begin{equation}
\label{weq5}
	g(x)=h_{C}\cdot f(x) \mbox{~~~for all~}x\in C.
\end{equation}
The element $h_{C}$ depends solely on $C$.

There is another notion of equality that is attractive but flawed.
Let $q:N\longrightarrow H\backslash N$ be the topological
quotient.  It is tempting to identify the morphisms of $f$ and $g$
if $q\circ f=q\circ g$ in the topological sense.  That is,
\begin{onecond}
\label{weq6}
\item	for each $x\in U$, there is $h_{x}\in H$ for which
	$g(x)=h_{x}\cdot f(x).$
\end{onecond}
Unlike (\ref{weq5}), the member $h_{x}\in H$ depends on each
point $x$.

If $q:N\longrightarrow H\backslash N$ is pseudo\'etale, then
conditions (\ref{weq5}) and (\ref{weq6}) are equivalent.  The
link is that $N_{0}[1,1]$ really is $q\times q$.  Thus,
$q\circ f=q\circ g$ implies existence of of $\delta =f\times g$
which is the basis of (\ref{weq5}).

Let us return to the explicit example of this Section.  Put
$f(v)=F(v)$ and $g(v)=F(-1\cdot v)$.  Then $F$ factors to a
morphism
$G\backslash M\longrightarrow H\backslash N$ if and only if
$f$ and $g$ determine the same morphism $M\longrightarrow N_{0}$.
Let $q:\bbr ^{2}\longrightarrow \{ -1,1\}\backslash \bbr ^{2}$
be the topological quotient.  Certainly $q\circ f=q\circ g$
in \cattop .  But, by inspection, $f$ and $g$ fail (\ref{weq5}).
In fact, $q$ is a false friend, and $F$ does not properly factor.

Equality should be a local property.  Let $\cal U$ be a family
of open subsets of $\bbr ^{3}$.  If, for each $V\in {\cal U}$,
restrictions of $f$ and $g$ to V are equal, then the restrictions
of $f$ and $g$ to $U=\cup {\cal U}$ should agree.  In fact, this is true.
However, one must realize that $U$ excludes lots of points.

Let $x\in \bbr ^{3}$.  Suppose that for each neighborhood $V$ of
$x$ there are $v,w\in V$ such that
$$\begin{array}{l}
	F(v)=F(-1\cdot v)\neq (-1)\cdot F(-1\cdot v), \\
	F(w)=(-1)\cdot F(-1\cdot w)\neq F(-1\cdot w).
	\end{array}$$
Then restrictions of $f$ and $g$ to any neighborhood of $x$
will not agree.  Obviously, every $x$ or the form $(0,0,z)$ has
this property; the function $F$ was tailored to be pathological
near such points.  Less obvious is the fact that each $x\in C_{n}$
has the same eccentricity.  In fact, the troublesome points divide
$\bbr ^{3}$ into disconnected chunks.  Indeed, if they did not, the
morphism $F$ would have factored.

\end{document}